 \documentclass[final]{IEEEtran}

\IEEEoverridecommandlockouts              

\usepackage{graphicx}

\usepackage{url}
\usepackage{amsmath,amssymb}
\usepackage{amsthm,thmtools}
\usepackage{amsfonts} 

\usepackage{enumerate}
\usepackage{color}                  

\usepackage{verbatim}

\usepackage{amssymb}
\usepackage{amsbsy}
\usepackage{amsmath}
\usepackage{setspace}
\usepackage{url}

\usepackage{float}
 \usepackage{dsfont}

 \newcommand{\dist}{\operatorname{dist}}

 \newcommand{\real}{\operatorname{Re}}
 
 \newcommand{\diag}{\operatorname{diag}}

\newcommand*\diff{\mathop{}\!\mathrm{d}}

\newcommand{\trace}{\operatorname{trace}}

\declaretheorem[name={Example},qed={\lower-0.3ex\hbox{$\square$}} ] {Example}

\declaretheorem[name={Definition}  ] {Definition}
\declaretheorem[name={Theorem}  ] {Theorem}
\declaretheorem[name={Lemma}  ] {Lemma}
\declaretheorem[name={Remark}  ] {Remark}
\declaretheorem[name={Corollary}  ] {Corollary}

\declaretheorem[name={Proposition}  ] {Proposition}

\newcommand {\R}{\mathbb R}

\newcommand {\C}{\mathbb C}

\newcommand{\be}{\begin{equation}}
\newcommand{\ee}{\end{equation}}

 \newcommand{\Z}{\mathcal Z}

\usepackage{lineno}

\begin{document}

%

%
\title{Generalization of the multiplicative and additive compounds of square matrices
and contraction in the Hausdorff dimension 
\thanks{The research of MM is
		supported in part 
		by a research grant  from  the Israel Science Foundation.  
}}

\author{Chengshuai Wu,  Raz Pines,  Michael Margaliot, and
 Jean-Jacques Slotine \thanks{
	\IEEEcompsocthanksitem
	 C. Wu and R. Pines are   with 
		the School of Electrical Engineering,
		Tel Aviv University, Tel-Aviv 69978, Israel.
		\IEEEcompsocthanksitem
		M. Margaliot (Corresponding Author) is  with the School of Electrical  Engineering,
		and the Sagol School of Neuroscience, 
		Tel-Aviv University, Tel-Aviv~69978, Israel.
		E-mail: \texttt{michaelm@tauex.tau.ac.il}
      \IEEEcompsocthanksitem
        J.-J. Slotine is with the Department of Mechanical Engineering and
        the Department of Brain and Cognitive Sciences, Massachusetts Institute of Technology, Cambridge, Massachusetts, USA.
        E-mail: \texttt{jjs@mit.edu}
        }}

\maketitle

\begin{center} 
		\today \\
\end{center}

\begin{abstract}
The $k$~multiplicative and $k$~additive compounds  of a matrix play an important role in geometry, multi-linear algebra,   the asymptotic analysis of nonlinear  dynamical systems, and in bounding   the Hausdorff dimension  of fractal sets. These compounds are defined for integer values of~$k$. Here, we introduce generalizations called the $\alpha$~multiplicative and $\alpha$~additive compounds of a square matrix, with~$\alpha$ real. We study the properties of these new compounds and demonstrate an application in the context of the   Douady and  Oesterl\'{e} Theorem. This leads to a generalization of  contracting  systems to $\alpha$ contracting  systems, with~$\alpha$ real. Roughly speaking, the dynamics of such systems contracts   any set with Hausdorff dimension larger than~$\alpha$. For~$\alpha=1$ they reduce to standard contracting systems. 
\end{abstract}


\section{Introduction}

Consider a matrix~$A \in \C^{n\times m }$, and  fix
an integer~$k \in \{1,\dots,\min\{m,n\}\}$.
 The~$k$ multiplicative compound matrix of~$A$,
 denoted~$A^{(k)}$, is the~$\binom{n}{k} \times \binom{m}{k}$ matrix that  includes all the minors of order~$k$ of~$A$ organized in a  lexicographic order. For example, if~$n=3$,~$m=2$, and~$k=2$, then~$A^{(2)}\in \R^{3\times 1}$
  and is given by 
 \[
 A^{(2)}=\begin{bmatrix} 
 A(\{1,2\} | \{1,2\} ) \\
 A(\{1,3\} | \{1,2\} )\\
 A(\{2,3\} | \{1,2\} )
 \end{bmatrix},
 \]
 where~$A(\alpha|\beta)$ denotes the minor of~$A$
 obtained by taking the rows indexed by~$\alpha$
 and the columns indexed by~$\beta$. 
Note that this implies that~$A^{(1)}=A$, and if~$m=n$ then~$A^{(n)}=\det(A)$.

 The Cauchy-Binet formula~\cite{total_book} asserts that for any $A \in \C^{n \times m}$, $B \in \C^{m \times r}$, and any~$k \in \{ 1,\dots ,\min\{n, m, r \}\}$,
 \be \label{eq:mulk}
 (AB)^{(k)} = A^{(k)} B^{(k)}.
 \ee
 This justifies the term multiplicative compound.  For example, if~$A,B$ are~$n\times n$
 then~\eqref{eq:mulk} with~$k=n$ reduces to the familiar formula~$\det(AB)=\det(A)\det(B)$.

 When~$n=m$, i.e.~$A$ is a square matrix, the $k$ additive compound matrix  of~$A$ is defined by
 \be \label{eq:defadd}
 A^{[k]}:= \frac{d}{d\varepsilon}
 (I+\varepsilon A)^{(k)} |_{\varepsilon=0}.
 \ee
 This implies that
 \be \label{eq:addtay}
 (I+\varepsilon A)^{(k)}=I+\varepsilon A^{[k]} +o(\varepsilon),
 \ee
i.e. $A^{[k]}$ is the first-order  term  in the Taylor series of~$(I+\varepsilon A)^{(k)}$. For example, if~$A\in\C^{n\times n}$ is diagonal, denoted~$A=\diag(\lambda_1,\dots,\lambda_n)$, and~$k=2$ then 
\begin{align*}
 (I & +\varepsilon A)^{(2)} \\
  & = (\diag(1+\varepsilon \lambda_1,\dots,1+\varepsilon \lambda_n  )   )^{(2)}\\
 &=\diag(  (1+\varepsilon \lambda_1)   (1+\varepsilon \lambda_2) ,
 \dots,(1+\varepsilon \lambda_{n-1})   (1+\varepsilon \lambda_n)
 ) \\
 &=I+\varepsilon \diag(\lambda_1+\lambda_2,\dots,\lambda_{n-1}+\lambda_n) +o(\varepsilon),
\end{align*}
so
\[
A^{[2]}= \diag(\lambda_1+\lambda_2,\dots,\lambda_{n-1}+\lambda_n).
\]
Note that every eigenvalue of~$A^{[2]}$ is the sum of two eigenvalues of~$A$.

It can be shown~\cite{muldo1990} that~\eqref{eq:mulk} and \eqref{eq:addtay}
imply that  
\be \label{eq:akadd}
(A + B)^{[k]} =  A^{[k]} + B^{[k]},
\ee
for any~$A, B \in \C^{n \times n}$.
This justifies the term additive compound. 

Compound matrices have found numerous applications in multi-linear algebra, geometry, and dynamical systems theory.
We quickly review some examples. 

The exterior product (or wedge product) of vectors
generalizes the notions of unsigned area  [volume] in dimension~$2$ [$3$] to an arbitrary dimension~$k$~\cite{wini2010}. 
For~$k $ vectors~$u^1,\dots,u^k\in\R^n$, with~$k\leq n$, the wedge product, denoted $u^1\wedge\dots\wedge u^k$ or~$\wedge_{i=1}^k u^i$, can be defined using the~$k$ multiplicative compound 
as  
\be \label{eq:defwedge}
u^i :=  \begin{bmatrix} u^1&u^2& \dots&u^k \end{bmatrix}^{(k)}.
\ee
Note that this has dimensions~$\binom{n}{k} \times \binom{k}{k}=\binom{n}{k} \times 1
$, i.e. it is a vector in~$\R^ {\binom{n}{k}} $. The magnitude of $u^1\wedge\dots\wedge u^k$, i.e., $\vert u^1\wedge\dots\wedge u^k \vert$, can be interpreted as the hyper volume of the $k$-dimensional parallelotope with edges~$u^1, \dots, u^k$.
For example, if~$a,b\in \R^3$ then
\begin{align*}
a\wedge  b = \begin{bmatrix}
a& b
\end{bmatrix}^{(2)} 
= \begin{bmatrix}
 a_1 & b_1 \\ a_2 & b_2 \\ a_3 & b_3
\end{bmatrix}^{(2)} = 
\begin{bmatrix}
 a_1 b_2 - b_1 a_2 \\
 a_1 b_3 - b_1 a_3 \\
 a_2 b_3- b_2 a_3 
\end{bmatrix},
\end{align*}
and the entries here are (up to a sign change) the entries of the cross product~$a \times b$. 
Recall that the magnitude of~$a\times b$ can be interpreted as the positive area of the parallelogram having~$a$ and~$b$ as sides.

Recall that~$A\in \R^{n \times m}$ is called   totally non-negative [totally positive]  if all its minors are non-negative [positive]. These matrices have found numerous applications~\cite{total_book,pinkus}. 
Clearly,~$A\in \R^{n \times m}$ is     totally non-negative [totally positive]  if and only if  
the multiplicative compound  matrices~$A^{(1)},A^{(2)},\dots,A^{(\min\{n,m\})}$ all have non-negative  [positive] entries. Thus,  every~$A^{(i)}$
can be analyzed using  the Perron-Frobenius theory of matrices with non-negative entries~\cite{matrx_ana}. This 
simple fact has important applications in the analysis of totally non-negative matrices (see, e.g.,~\cite{gk_book,rola_spect}).

We now describe some applications of compound matrices in dynamical systems theory.
In this context,
the relevant case is square matrices.
Suppose that~$X(t)$ is the solution of the linear matrix differential equation
\be\label{eq:xax}
\frac{d}{dt}  X(t) =A(t) X(t),\quad X(t_0)=X_0,
\ee
with~$A:\R  \to \R^{n \times n }$ continuous. Then 
\be\label{eq:xk}
\frac{d}{dt} X^{(k)}(t) = A^{[k]}(t)  X^{(k)}(t),\quad X^{(k)}(t_0)=(X_0)^{(k)} ,
\ee
for any~$k\in\{1,\dots,n\}$.
In other words, all the~$k$ minors of~$X$, stacked in the matrix~$X^{(k)}$, also follow a linear dynamics, with the matrix~$A^{[k]}$ (see, e.g.,~\cite{muldo1990}). 
Note that if~$A$ is time-invariant and~$t_0=0$ then the solution of~\eqref{eq:xax}
is~$X(t)=\exp(At)X_0$
so~$(X(t))^{(k)}=(\exp(At))^{(k)} (X_0)^{(k)}$ and combining this with~\eqref{eq:xk}
gives
\be\label{eq:expat}
(\exp(At))^{(k)}= 
\exp(A^{[k]} t),\text { for all } t\in\R. 
\ee
Eq.~\eqref{eq:xk} has important applications  in the analysis of time-varying nonlinear dynamical systems in the form
\be\label{eq:nonlin}
\dot x(t)=f(t,x(t)),
\ee
where~$f:\R\times \R^n\to\R^n$ is~$C^1$.
Indeed, let~$x(t,a)$ denote the solution at time~$t\geq 0 $ of~\eqref{eq:nonlin} with~$x(0)=a$, and let
\be \label{eq:jacf}
J_f(t,x):=\frac{\partial }{\partial x} f(t,x)
\ee
denote the Jacobian of the vector field~$f$. Then the variational equation associated with~\eqref{eq:nonlin}
along the trajectory~$x(t,a)$ is
\be\label{eq:var}
\dot y(t)=J_f(t,x(t,a)) y(t).
\ee
Analysis of this linear time-varying  equation plays an important role in the asymptotic analysis of~\eqref{eq:nonlin}. 
Combining this with~\eqref{eq:xk} has   far-reaching applications in the theory of nonlinear dynamical systems~\cite{muldo1990}.
Recent applications include:
\begin{itemize}
    \item 
totally positive differential systems~\cite{schwarz1970,margaliot2019revisiting}, that is,
systems where the
transition matrix corresponding to the variational  equation~\eqref{eq:var}
is~TP (see also~\cite{CTPDS});
\item
 $k$-cooperative dynamical systems~\cite{Eyal_k_posi}, that is,
 systems where~$J_f^{[k]}$ is a Metzler matrix;  
 \item
$k$-order contracting systems~\cite{kordercont}, that is,
systems where~$J_f^{[k]}$ is infinitesimally contracting;
\item
the notion of a discrete-time $k$-diagonally stable dynamical system, that is, a system of the form~$x(k+1)=Ax(k)$,  and there exists a positive-definite  diagonal matrix~$D$,  such that
$
(A^{(k)} )^ T D A^{(k)} -D
$
is negative-definite~\cite{cheng_diag_stab}.
\end{itemize}

Since the $k$ multiplicative compound is based on collecting all the $k\times k$ minors of a matrix, it is naturally defined only for integer values of~$k$.

Here, we introduce a generalization of the~$k$  multiplicative and~$k$ additive
compound  of a square  matrix, called the $\alpha$ multiplicative compound
and   $\alpha$ additive compound,
where~$\alpha\geq 1$ is a \emph{real} number.
For~$k<\alpha<k+1$   the $\alpha$ compounds may be interpreted as a weighted   interpolation of the $k$ and~$k+1$ compounds. 
When~$\alpha$ is an integer this (almost) reduces to the standard~$k$ compounds. 
This generalization  is motivated by the 
Hausdorff dimension of a set
and, in particular, the seminal 
 Douady and Oesterl\'{e} Theorem~\cite{Douady1980} that  provides an upper bound for the Hausdorff dimension of a set that is negatively invariant under a $C^1$ mapping.
As an application, we show that the~$\alpha$ compounds can be used to  provide elegant and intuitive expressions for the basic terms that appear in 
this theorem. Furthermore, this naturally  leads to the new notion of~$\alpha$ contracting systems, with~$\alpha$ real, which generalizes the notion of~$k$-order contracting systems with~$k$ an integer~\cite{leonov,weak_slotine,kordercont}. We analyze the properties of $\alpha$~contracting systems and demonstrate their applications. 
Our results 
show that if an $n$-dimensional  dynamical system contracts $n$-dimensional volumes then there exists a minimal \emph{real} value~$\alpha^* \in[1,n]$ such that the system is~$\alpha$ contracting  for any~$\alpha>\alpha^*$. 
Roughly speaking, an~$\alpha^*$ contracting system  contracts any set with a Hausdorff dimension larger than~$\alpha^*$.
This generates (in a given metric) a continuum of contraction  
 instead of the standard binary view, namely, that a system is either contracting or not contracting.

The remainder of this paper is organized as follows. The next section reviews several known definitions and   results that will be used later on. Section~\ref{sec:main} describes our main results. These include the definitions of the~$\alpha$ multiplicative compound and~$\alpha$
additive compound 
of a matrix, and analysis of the  properties of these compounds. 
Section~\ref{sec:app} describes an application of these 
new notions and introduces $\alpha$~contracting systems.
The last section concludes and describes several directions for future research.
We focus here on the theory of the generalized compounds and $\alpha$ contracting systems, leaving applications to a sequel paper. 

\section{Preliminaries}
To make this paper more self-contained, we 
briefly review several topics that are  needed
to define,  analyze and apply  the $\alpha$ compounds, and~$\alpha$ contracting systems. 
We  begin with   reviewing the Hausdorff dimension 
of a set, following~\cite{smith_hauss,Pogromsky_2000,book_volker2021}.

\subsection{Hausdorff dimension  }

Let~$K$ be a  set in~$\R^n$.
For~$d,\epsilon>0$, the $d$-measured volume of~$\epsilon$-coverings of~$K$ is:
\begin{align*}
\zeta(K,\epsilon,d):=  \inf & \{
\sum_i r_i^d : \text{there exists a countable} \\ & \text{cover of }  K \text{ by balls with radii } r_i \leq \epsilon \}.
\end{align*}
Note that the covering may include balls of different sizes, but all are bounded by~$\epsilon$.
Note also that if~$K$ is compact  then it would suffice to use finite coverings, since every open cover of~$K$ has a finite subcover.

By definition,~$\zeta(K,\epsilon,d)$ is non-increasing in~$\epsilon$. 
The Hausdorff $d$-measure of~$K$ is
\[
m(K,d):=\lim_{\epsilon \downarrow 0}\zeta(K,\epsilon,d),
\]
where the limit may be infinite.

For any~$s>0$, we have
\[
\zeta(K,\epsilon,d+s)\leq \epsilon^s
\zeta(K,\epsilon,d) ,
\]
implying that if~$m(K,d) <\infty$
then~$m(K,d) =0$. Thus, there is a unique~$d^*   \in [0,n]$ such that
$m(K,d)=0$ for all~$d>d^*$,
and~$m(K,d)=\infty$ for all~$d<d^*$.
The  Hausdorff dimension of~$K$
is 
\[ 
\dim_H K:=d^* .
\]
Intuitively,  if we try to cover a square (which is a 2D set) by 1D balls (i.e. lines) lines  then we need an infinite number of lines, but once we try a cover with 2D balls, the number of balls needed is finite.
So,~$d^*$ is exactly
the
  dimension for which 
the ``volume'' of~$K$ becomes finite.

For smooth  shapes, or shapes with a small number of ``corners'',  the Hausdorff dimension is an integer agreeing with the more standard  topological dimension.
For example, suppose that~$K$ is an $\ell$-dimensional cube  in~$\R^n$. Intuitively speaking, for any~$\epsilon>0$
we require~$\theta\left ((1/\epsilon  )^\ell\right )$
 balls of radius~$\epsilon$
 to cover~$K$. 
Hence,
\begin{align}\label{eq:nbhpw}
\zeta(K,\epsilon,d)\approx  (1/\epsilon)^\ell \epsilon^d  .
\end{align}
As~$\epsilon\downarrow 0$, the right-hand side of~\eqref{eq:nbhpw} goes to~$\infty$ if~$d<\ell$, and to zero if~$d>\ell$. 
It is not difficult to see that using balls with varying sizes does not change the analysis, so
\[
\dim_H K=\ell. 
\]

However, 
 for  fractal sets (e.g. sets that contain strange  
 attractors of chaotic  dynamical  systems) the Hausdorff dimension is typically not an  integer.\footnote{In fact, one possible definition  of a fractal set  is:
a set whose  Hausdorff dimension strictly
exceeds its topological dimension~\cite{Haus_dim_general}.}
In this context,
the   Hausdorff dimension is useful in
quantifying sets of Lebesgue measure zero which are nevertheless ``substantial''.
The next example from~\cite{cantor93}
demonstrates this.
\begin{Example}
 The Cantor set~$E\subset[0,1]$ is defined inductively as follows.
 Let~$E_0:=[0,1]$. For~$j\geq 1$, the set~$E_j$ is obtained  by removing the open middle third in any interval 
 in~$E_{j-1}$. For example,
 $E_1=[0,1/3]\cup[2/3,1]$, and
 $E_2=[0,1/9]\cup[2/9,1/3]\cup[2/3,7/9]\cup[8/9,1] $. Cantor's set is~$E:=\cap_{j=0}^\infty E_j$. Each~$E_j$ is the union of~$2^j$ intervals of length~$3^{-j}$.
 The topological dimension 
 of~$E$  is thus
 \[
 \lim_{j\to\infty} (2/3)^j=0.
 \]
 We now determine~$\dim_H E$. 
 It can be shown~\cite[Ch.~3]{book_volker2021}  that it is enough to consider the  cover of~$E_j$ by~$2^j$ intervals, each  of length~$3^{-j}$, so for any~$\epsilon>0$ sufficiently  small, $ \zeta(E_j,\epsilon,d)=2^j 3^{-jd}$.
 Thus,~$m(E_j, d)=\left( \frac{2}{3^d}\right)^j$
 and
 \[
m(E , d)=\lim_{j\to\infty }\left( \frac{2}{3^d}\right)^j.
 \]
 If~$d>\log(2)/\log(3)$ then~$m(E,d)=0$, and if~$d<\log(2)/\log(3)$ then~$m(E,d)=\infty$.
 Thus,
 $
 \dim_H  E =\log(2)/\log(3)\approx 0.631.
 $
 Intuitively, this implies that the Cantor set 
is less than a line, but more than a discrete set of points.
\end{Example}
 
 The next result summarizes some useful  properties of $\dim_H$.
 \begin{Proposition}[{\cite[Chapter~3]{book_volker2021}}]
 The Hausdorff dimension satisfies the following properties:
 \begin{itemize}
 \item~$\dim_H \emptyset =0 $;
     \item monotonicity:  If~$A,B\subseteq \R^n$ with~$ A\subseteq B$ then~$\dim_H A \leq \dim_H B$;
     \item countable subadditivity: If~$A_i \subseteq \R^n$, $i=1,2,\dots$,  then
     $
     \dim_H (\cup_i A_i)\leq \sum_i \dim_H (A_i);
    $
    \item If~$A,B\subseteq \R^n$ are such that
    \[
    \inf\{\dist(x,y):x\in A, y\in B\} >0
    \]
    then
    \[
    \dim_H(A\cup B)=\dim_H A+\dim_H B.
    \]
 \end{itemize}
 
 \end{Proposition}
 The first three properties imply that~$\dim_H$ is an outer measure, and the fourth one that it is a metric outer measure.

\subsection{Explicit formula for~$A^{[k]}$}
The additive compound,
  defined in~\eqref{eq:defadd},     can be given explicitly in terms of the entries~$a_{ij}$ of~$A$. 
\begin{Proposition} [\cite{schwarz1970,fiedler_book}] \label{prop:explicitak}
Let~$A\in\R^{n\times n}$ and fix~$k\in\{1,\dots,n\}$, 
$1\leq i_1<\dots <i_k \leq n$, and~$1\leq j_1<\dots< j_k \leq n$.
The entry of~$A^{[k]}$ corresponding to~$(\alpha|\beta)=(i_1,\dots,i_k|j_1,\dots,j_k) $  is:
\begin{itemize}
\item $\sum_{\ell=1}^k a_{i_\ell i_\ell}$ if~$i_\ell=j_\ell$ for all~$\ell\in\{1,\dots,k\}$;
\item $(-1)^{\ell+m} a_{i_\ell j_m}$ 	 if all the indices in~$\alpha$ and~$\beta$ coincide except for a single index~$i_\ell\not = j_m $; and
\item $0$, otherwise. 
\end{itemize}
\end{Proposition} 

To explain this, consider for example the case~$k=2$. Let~$B:=(I+\varepsilon A)^{(2)}$, and let
\[
\delta_{pq}:=\begin{cases}1, \text{ if } p=q,\\
0, \text{ otherwise.}
\end{cases}
\]
Fix~$1\leq i_1<i_2\leq n$  and~$1\leq j_1<j_2\leq n$.
Then
\begin{align*}
    B(i_1,i_2|j_1,j_2) = &b_{i_1j_1}b_{i_2j_2}-b_{i_1j_2}b_{i_2j_1}\\
    =&(\delta_{i_1j_1}+ \varepsilon a_{i_1j_1} )
     (\delta_{i_2j_2}+ \varepsilon a_{i_2 j_2} ) \\
     &-(\delta_{i_1j_2}+ \varepsilon a_{i_1j_2} )
     (\delta_{i_2j_1}+ \varepsilon a_{i_2 j_1} )\\
     =&c+   (\delta_{i_1j_1} a_{i_2 j_2} + \delta_{i_2 j_2} a_{i_1 j_1} -\delta_{i_1 j_2} a_{i_2 j_1}  \\&- \delta_{i_2 j_1} a_{i_1 j_2})\varepsilon 
     +o(\varepsilon),
\end{align*}
where~$c$ is a constant that does not depend on~$\varepsilon$. Thus,~\eqref{eq:defadd} gives
\begin{align*}
 A^{[2]}(i_1,i_2|j_1,j_2
 )=& \delta_{i_1j_1} a_{i_2 j_2} + \delta_{i_2 j_2} a_{i_1 j_1} \\& -\delta_{i_1 j_2} a_{i_2 j_1} - \delta_{i_2 j_1} a_{i_1 j_2},
 \end{align*}
and it is straightforward to see that this agrees with the expression given in 
Prop.~\ref{prop:explicitak}.
Note that Prop.~\ref{prop:explicitak} implies, in particular,  that
$ A^{[1]}=A$,
and $A^{[n]} =\trace(A)$.

\subsection{Real power of a square non-singular matrix}
We first recall  the definition of  the real power of a complex number. Any complex number $a \in \C\setminus\{0\}$ can be written in the polar representation 
$
a = \vert a \vert  \exp(j \theta(a)),
$
where~$j:=\sqrt{-1}$,
$\vert a \vert >0$ is the  modulus of~$a$, and~$\theta(a) \in (-\pi , \pi]$ is the  argument 
of~$a$.
  Then for any~$\alpha \in \R $,  
\begin{equation}
a^\alpha :=  \vert a \vert^\alpha  \exp(j \alpha \theta(a)).
\end{equation}
For example, for~$a=-5$
we have
$
(-5)^\alpha=5^\alpha \exp( j\alpha \pi )$.
Note that although~$-5$ is real, $(-5)^\alpha$ is in general a complex (non-real) number.

Recall that any $A \in \C^{n \times n}$ admits a Jordan canonical form~\cite{achieser1992theory}: there exist $T, J \in \C^{n \times n}$, with~$T$ non-singular, 
such that
\be \label{eq:jordan}
A =  T^{-1} J T, \quad J = \diag(J_1, J_2, \dots, J_p),
\ee  
where every~$J_i$ is a Jordan block in the   form
\begin{equation}  \label{eq:jordanblock}
J_i = 
\begin{bmatrix}
\ell_i & 1 & ~& \\
~ & \ell_i & \ddots & ~ \\
~& ~ & \ddots & 1 \\
~ & ~& ~& \ell_i
\end{bmatrix} \in \C^{m_i \times m_i},
\end{equation}
with~$ \sum_{i=1}^p m_i = n$, and every~$\ell_i$, $i=1,\dots,p$,
is an eigenvalue of~$A$.
The matrix~$J$ is unique, up    to the ordering of the blocks~$J_i$.

Since the real power of square matrices is a particular class of a  matrix function  \cite{gantmacher1966theory, higham2008functions}, it is defined according to the general definition given in~\cite[Def. 1.2]{higham2008functions}.

\begin{Definition}[Real power of a non-singular square  matrix] \label{def:rpsm}
Consider a non-singular matrix $A \in \C^{n \times n}$, given in the Jordan canonical form~\eqref{eq:jordan}, and let~$\alpha \in \R$. Then, 
\be
A^\alpha := T^{-1} J^\alpha T,
\ee
where
\be
J^\alpha := \diag(J_1^\alpha, \dots, J_p^\alpha  ),  
\ee
with
\be \label{eq:jalp}
J_i^\alpha :=
\begin{bmatrix}
\ell_i^\alpha &  \frac{\alpha \ell_i^{\alpha-1}}{1!} & \cdots & \frac{\prod_{j=0}^{m_i-2} (\alpha-j) \ell_i^{\alpha-m_i +1}  }  {(m_i -1) !} \\
~ & \ell_i^\alpha & \ddots & \vdots \\
~ & ~& \ddots & \frac{\alpha \ell_i^{\alpha-1}}{1!}  \\
~ & ~& ~& \ell_i^\alpha
\end{bmatrix}.
\ee
\end{Definition}

The next result describes some of the properties of~$A^\alpha$.
\begin{Proposition}[{\cite[Thm.~1.13, Thm~1.15]{higham2008functions}}]\label{prop:real}
Let~$A \in\C^{n\times n}$ be non-singular
    with the  Jordan   canonical form~\eqref{eq:jordan}. Fix~$\alpha,\beta\in \R$. Then
    \begin{enumerate}[(a)]
   \item the    eigenvalues of $A^\alpha$ are $\ell_i^\alpha$, $i = 1, \dots, p$;\label{pp:realeig}  
    \item \label{pp:realT} $(A^T)^\alpha = (A^\alpha)^T$; 
    \item \label{pp:realinv} $(X A X^{-1})^\alpha = X A^\alpha X^{-1}$, for any non-singular matrix~$X\in\C^{n\times n}$; 
    \item   $A^\alpha A^\beta = A^{\alpha + \beta}$. 
    \end{enumerate}
\end{Proposition}

Note that $A^\alpha$ is not necessarily real even if $A$ is real. The next result provides a sufficient condition guaranteeing  that~$A^\alpha$ is real for any~$\alpha\in\R$. 
Let $\R_{\leq 0} : =\{x \in \R ~| ~ x \leq 0\}$, and~$\R_{> 0} : =\{x \in \R ~| ~ x > 0\}$.
 Define the set of matrices~$\Omega_n : = \{ X \in \R^{n \times n} ~|~ \text{spec}(X) \cap \R_{\leq 0} = \emptyset \}$, where $\text{spec}(X)$  is  the set of eigenvalues of $X$.
\begin{Proposition} \label{prop:sufreal}
  If $A \in \Omega_n$, then $A^\alpha \in\R^{n\times n}$   for any $\alpha \in \R$.
\end{Proposition}

\begin{IEEEproof}
For~$a\in \C$, let~$\overline a$ denote the complex conjugate of~$a$. 
 Note that $A \in \Omega_n$ implies that $A$ is non-singular. Definition~\ref{def:rpsm} guarantees that $A^\alpha$ is well-defined. Furthermore, the function $f(a) := a^\alpha$ is analytic on $\C\setminus \R_{\leq 0}$, and any connected component of~$\C\setminus \R_{\leq 0}$ is closed under conjugation. Note that $f(\R \cap (\C\setminus \R_{\leq 0})) = f(\R_{>0}) \subset \R$. Therefore, the conditions given  in \cite[Thm. 1.18]{higham2008functions} are satisfied. This ensures that $f(\overline{A}) =  \overline{f(A)}$. Since $A$ is real, i.e., $\overline{A} = A$, this implies that $f(A) = \overline{f(A)}$, that is, $f(A)$ is  real.
\end{IEEEproof}

The next result provides a formula for 
the derivative of   a  
real power of a parameter-dependent    matrix.
This result will be used later on to 
derive a simple expression for the generalized additive compound of a matrix. 

\begin{Lemma}\label{lem:dotAalp}
Consider a  matrix-valued mapping $A: \R \to \C^{n \times n}$. Assume that $A(\varepsilon)$ is non-singular and has constant (generalized) eigenvectors, that is, it can be written in the Jordan canonical form 
\[
A(\varepsilon) = T^{-1} J(\varepsilon)T,~ J(\varepsilon) = \diag(J_1(\varepsilon), J_2(\varepsilon), \dots, J_p(\varepsilon)),
\] 
where $T \in \C^{n \times n}$, and every~$J_i: \R \to \C^{m_i  \times m_i }$,~$\sum_{i=1}^p m_i = n$, is a Jordan block as described in~\eqref{eq:jordanblock}, where $\ell_i: \R \to \C\setminus\{0\}$, $i = 1, \dots, p$, is $C^1$.
Then for any $\alpha \in \R$,
\be\label{eq:eqtwen}
\frac{\diff }{\diff \varepsilon} A^{\alpha}(\varepsilon) = \alpha A^{\alpha -1}(\varepsilon) \frac{\diff }{\diff \varepsilon} A(\varepsilon) .
\ee
\end{Lemma}

\begin{IEEEproof}
By \eqref{eq:jalp},
$
\frac{\diff }{\diff \varepsilon} J_i^\alpha(\varepsilon) = \alpha J_i^{\alpha -1 } (\varepsilon)\frac{\diff }{\diff \varepsilon} \ell_i(\varepsilon) 
$, so 
\begin{align*}
\frac{\diff }{\diff \varepsilon} J^\alpha(\varepsilon) =& \alpha J^{\alpha -1 }(\varepsilon) \diag \left  (  \frac{\diff }{\diff \varepsilon} \ell_1(\varepsilon) I_{m_1}, \dots, \frac{\diff } {\diff \varepsilon} \ell_p(\varepsilon) I_{m_p} \right  ) \\
=&  \alpha J^{\alpha -1 }(\varepsilon)  \frac{\diff }{\diff \varepsilon} J(\varepsilon),
\end{align*}
where $I_{s}$ denotes the $ s \times s$ identity matrix. Thus, 
\begin{equation*}\begin{split}
\frac{\diff }{\diff \varepsilon} A^{\alpha}(\varepsilon) =& T^{-1}  \frac{\diff }{\diff \varepsilon} J^\alpha(\varepsilon) T \\
=&  \alpha T^{-1}   J^{\alpha -1 }(\varepsilon)  T T^{-1} \frac{\diff }{\diff \varepsilon}
J(\varepsilon) T,  
\end{split}
\end{equation*}
and using  Definition~\ref{def:rpsm} yields~\eqref{eq:eqtwen}.
\end{IEEEproof}
 
 \subsection{Kronecker product and Kronecker sum of matrices} 
 The Kronecker product
 of 
two matrices~$A \in \C^{n \times m}$ and~$B \in \C^{p \times q}$ is
\be
A \otimes B := 
\begin{bmatrix}
a_{11} B & a_{12}B & \cdots & a_{1m} B  \\
a_{21} B & a_{22}B & \cdots & a_{2m} B \\
\vdots & \vdots  & \ddots & \vdots \\
a_{n1} B & a_{n2}B & \cdots& a_{nm} B
\end{bmatrix},
\ee
where $a_{ij}$ denotes the $ij$th entry of $A$. Hence, $A \otimes B \in \C^{(np) \times (mq)}$. 

The Kronecker sum of two
square matrices~$X \in \C^{n \times n}$ and~$Y \in \C^{m \times m}$ is
 \begin{equation}
    X \oplus Y :=  X \otimes I_m + I_n \otimes Y,
 \end{equation}
 where $I_r$ denotes the $r \times r$ identity matrix. 
 
 We list several  properties of the Kronecker product and Kronecker sum that will be used later on.   
 For~$A\in \C^{n\times m}$, let~$A^*$ to
 denote  the conjugate transpose of $A$.
 
 \begin{Lemma}[see e.g. \cite{graham2018kronecker}]  \label{lem:pp}
Consider   matrices $A, C \in \C^{n \times m}$, $B, D \in \C^{p \times q}$, $F \in \C^{m \times \ell}$, $G \in \C^{q \times r}$,  $X \in \C^{n \times n}$, and~$Y \in  \C^{m \times m}$. Then, 
 \begin{enumerate}[(a)]
     \item \label{pp:asos} $(c A) \otimes B = A \otimes (c B) =  c (A \otimes B)$ for any $c \in \C$;
     \item \label{pp:dist} $(A + C) \otimes B  =  A \otimes B + C \otimes B$ ;
     \item \label{pp:dist1} $A \otimes (B + D) = A \otimes B + A \otimes D$;
     \item \label{pp:mixpro} $(A \otimes B) (F \otimes G) = (AF)\otimes (BG)$;
     \item  $(A \otimes B)^* = A^* \otimes B^* $;
     \item \label{pp:tranp} $(A \otimes B)^T = A^T \otimes B^T $;
    \item \label{pp:inverse} if $X,Y$ are non-singular then $(X \otimes Y)^{-1}   =  X^{-1} \otimes Y^{-1}$;
     \item \label{pp:eigotimes} Let $\lambda_i(X)$, $i = 1, \dots, n$, and $\lambda_j(Y)$, $j = 1, \dots, m$, denote the eigenvalues of $X$ and $Y$, respectively. Then, $X \otimes Y$ has eigenvalues  $\lambda_i(X)\lambda_j(Y)$,  $i = 1, \dots, n$, $j = 1, \dots, m$;
     \item \label{pp:eigoplus}  $X \oplus Y$ has eigenvalues $\lambda_i(X) + \lambda_j(Y)$,  $i = 1, \dots, n$, $j = 1, \dots, m$;
     \item \label{pp:expkro} $ \exp(X) \otimes \exp(Y) = \exp(X \oplus Y)$.
 \end{enumerate}
 \end{Lemma}
 
 Property \eqref{pp:asos} implies 
 that we can  write $c A \otimes B := (c A) \otimes B$ or $A \otimes (c B)$, without  any ambiguity.

For a real number $p \geq 1$, let 
\be \label{eq:deflpnorm}
 \vert x \vert_p : = (\vert x_1 \vert^p + \cdots +\vert x_n \vert^p)^{\frac{1}{p}},
\ee
denote   the $L_p$ vector norm of $x \in \C^n$,
and let~$||\cdot ||_p$  denote the   induced matrix norm. 
Recall~\cite{norm_kron_prod} that
a norm~$|\cdot|:\C^n\to\R_+$ is called monotonic if for any~$x,y\in\C^n$ such that
$
|x_i|\leq |y_i| $ for all $i\in\{1,\dots,n\},
$
we have~$|x|\leq |y|$.
The $L_p$ norm  for any $p \geq 1$ is monotonic.
  
The next result uses  ideas from~\cite{norm_kron_prod} to determine   the induced~$L_p$ matrix norms of a Kronecker product of two matrices.

\begin{Proposition} \label{prop:normkrop}
Fix~$A \in \C^{\ell \times n}$, $B \in \C^{r \times m}$,  and~$p\geq 1$. Then, 
\be
\Vert A \otimes B \Vert_p = \Vert A \Vert_p \Vert B \Vert_p.
\ee
\end{Proposition}
\begin{IEEEproof}
It is easy to verify that
 \be \label{eq:cross}
 \vert x \otimes y \vert_p = \vert x \vert_p \vert y \vert_p,  \text{ for all } x \in \C^n,\; y \in \C^m.
 \ee
   Hence, in the terminology of~\cite{norm_kron_prod},
the~$L_p$ norms are   \emph{cross norms}.
 Let 
 \[
 \C^n \otimes \C^m := \{  z = x \otimes y: x \in \C^n, y \in \C^m \}.
 \]
 Note that $\C^n \otimes \C^m \subseteq \C^{nm}$. Then, 
 \begin{align*}
 \Vert A \otimes B \Vert_p =& \sup_{z \in \C^{nm}} \frac{\vert (A \otimes B) z \vert_p}{\vert z \vert_p} \\
 \geq & \sup_{z \in \C^n \otimes \C^m} \frac{\vert (A \otimes B) z \vert_p}{\vert z \vert_p} \\
 =& \sup_{x \in \C^n, y \in \C^m} \frac{\vert (A \otimes B) (x \otimes y) \vert_p}{\vert x \otimes y \vert_p}, 
 \end{align*}
 and applying Property~\eqref{pp:mixpro} in Lemma~\ref{lem:pp} and~\eqref{eq:cross}  yields
 \begin{align} \label{eq:kroabg}
 \Vert A \otimes B \Vert_p 
 \geq  \sup_{x \in \C^n} \frac{\vert Ax \vert_p}{\vert x\vert_p}\sup_{y \in \C^m} \frac{\vert By \vert_p}{\vert y\vert_p}  
 = \Vert A \Vert_p \Vert B \Vert_p.
 \end{align}
 Thus, to complete the proof we need to show that~$\Vert A \otimes B \Vert_p\leq \Vert A \Vert_p \Vert B \Vert_p$.
 Note that $A \otimes B = (A \otimes I_r) (I_n \otimes B) $. Since   induced matrix norms are sub-multiplicative, 
 \be \label{eq:submul}
 \Vert A \otimes B \Vert_p \leq \Vert A \otimes I_r \Vert_p \Vert I_n \otimes B \Vert_p.
 \ee
 Let $\{e^i\}_{i=1}^n$, $\{f^j\}_{j=1}^m$ denote the   canonical basis of~$\C^n$  and~$\C^m$, respectively. 
 Any~$z \in \C^{nm}$ can be written as
 \be \label{eq:zeifj}
 z = \sum_{i=1}^n e^i \otimes y^i \text{ and  } z=\sum_{j=1}^m x^j \otimes f^j ,
 \ee
 for some vectors $y^i \in \C^m$, $x^j \in \C^n$. As in~\cite{norm_kron_prod}, consider  the   norms
 \begin{align} \label{eq:defnormef}
    \vert z \vert_{p, e} :=   \vert \sum_{i=1}^n \vert y^i \vert_p e^i  \vert_p, \quad
    \vert z \vert_{p, f} :=   \vert \sum_{j=1}^m \vert x^j \vert_p f^j  \vert_p.
 \end{align}
It is easy to show that 
\be \label{eq:zpef}
\vert z \vert_{p, e} = \vert z \vert_{p, f}= \vert z \vert_p,
\text{ for any } z \in \C^{nm}.
\ee
 Now,
 \begin{align*}  
    \vert (I_n \otimes B)z \vert_{p,e} &= 
       \vert 
       (I_n\otimes B)(\sum e^i\otimes y^i)
       \vert_{p,e}
         \nonumber 
    \\
    &=
     \vert  \sum e^i\otimes B y^i \vert_{p,e} \\
     &= \vert 
     \sum \vert B y^i \vert_p
     e^i\vert_p\nonumber\\
    &\leq \Vert B \Vert_p
    \vert \sum \vert   y^i \vert_p
     e^i\vert_p\nonumber\\
     &=\Vert B \Vert_p \vert z \vert_{p,e}.\nonumber
 \end{align*}
  Thus,
 $
 \Vert I_n \otimes B \Vert_{p,e} \leq \Vert B \Vert_p.
$
 On the other hand,~\eqref{eq:kroabg} and~\eqref{eq:zpef} imply that $ \Vert I_n \otimes B \Vert_{p,e} = \Vert I_n \otimes B \Vert_p \geq \Vert I_n \Vert_p \Vert B \Vert_p = \Vert B \Vert_p$. Hence, 
 $
 \Vert I_n \otimes B \Vert_p = \Vert B \Vert_p.
 $
A similar argument using 
the norm~$\vert \cdot \vert_{p,f}$ yields
  $
 \Vert A \otimes I_r \Vert_p = \Vert A \Vert_p
 $.
 Now~\eqref{eq:submul} yields 
 $
 \Vert A \otimes B \Vert_p \leq \Vert A \Vert_p \Vert B \Vert_p
 $.
\end{IEEEproof}

\subsection{Matrix measures}
Let~$|\cdot|:\C^n\to\R_+$ denote a vector norm. The induced matrix norm~$||\cdot||:\C^{n\times n}\to\R_+$ is~$||A||:=\max_{|x|=1}|A x |$, and the induced matrix measure~$\mu:\C^{n
\times n}\to\R$ is  
\be \label{eq:defmu}
\mu(A):=\lim_{\varepsilon \downarrow 0}
\frac{||I+\varepsilon A||-1}{\varepsilon}.
\ee
Matrix measures (also called logarithmic norms) play an important role in numerical linear algebra~\cite{Strom_loga_norm} and in systems and control theory~\cite{LOHMILLER1998683, sontag_cotraction_tutorial}. The reason for this is two-fold. First,  if there exists a matrix measure such that the matrix~$A(t)$ in~\eqref{eq:xax} satisfies
\be\label{eq:xhepr}
\mu(A(t))\leq \eta, \text{ for all } t\geq 0, 
\ee
then
\[
||X(t)||\leq \exp(\eta t) ||X(0)||, \text{ for all } t\geq 0.
\]
In particular, if~$\eta<0$ then this implies exponential convergence to zero with rate~$\eta$.
Second, let $\mu_p$ denote the matrix measure induced by the~$L_p$ vector norm, with $p\in\{1,2,\infty\}$.
Then for any~$A\in\C^{n\times n}$,
\begin{align}\label{eq:mmeas}
    \mu_1(A)&= \max_{j\in\{1,\dots,n\}} \Big( \real(a_{jj})+\sum_{ \genfrac {}{}{0pt}{}
    {i=1}   {i\not =j}}^n |a_{ij}| \Big )  ,\nonumber \\
    \mu_2(A)&=\lambda_{\max}(A_\text{sym} ) , \\
    \mu_\infty(A)&= \max_{i \in\{1,\dots,n\}} \Big ( \real(a_{ii}) +\sum_{  
    \genfrac {}{}{0pt}{} {i=1}  {i\not =j} }^n |a_{ij}| \Big ) ,\nonumber
\end{align}
where $A_\text{sym} := (A + A^*)/2$,
 and $\lambda_{\max}(B)$ denotes the largest eigenvalue of a Hermitian matrix~$B$.
Using these explicit  formulas, it is sometimes possible to easily verify that~\eqref{eq:xhepr} indeed holds,   without   computing~$X(t)$ itself.

Let $Q^{k,n}$ denote the 
sequence of $k$-tuples of distinct numbers from $\{1, \dots, n \}$, in lexicographic   order.
For example,
\[
Q^{2,3}=\{  \{ 1,2\}, \{ 1,3\}, \{2,3\} \}.
\]
Note that there are $\binom{n}{k}$ such $k$-tuples. 
Let~$Q^{k,n}_\ell$ denote the~$\ell$th tuple in~$Q^{k,n}$. For example,
$Q^{2,3}_2=  \{ 1,3\} .
$

For $A \in \C^{n \times n}$ and any $k \in \{1 , \dots, n\}$, the matrix measures for~$A^{[k]}$ are~\cite{muldo1990}:
\begin{align}\label{eq:matirxm_k}
\mu_1(A^{[k]}) = &
 \max_{ \{i_1,\dots, i_k\}\in Q^{k,n} }\Big ( \sum_{p = 1}^k \real(a_{i_p,i_p})  \nonumber \\
 &+ \sum_{\substack{j \notin \{i_1,\dots, i_k\} }}(|a_{j,i_1}| + \cdots + |a_{j,i_k}|)\Big ) ,\nonumber \\
\mu_2(A^{[k]}) =& \sum_{i=1}^k \lambda_i\left( A_\text{sym} \right) , \\
\mu_{\infty}(A^{[k]}) =&
\max_{\{i_1,\dots, i_k\}\in Q^{k,n}}\Big(  \sum_{p=1}^k \real(a_{i_p,i_p})   \nonumber \\
&+ \sum_{\substack{j \notin \{i_1,\dots,i_k\} }}(|a_{i_1,j}| + \cdots + |a_{i_k,j}|)\Big).
\nonumber
\end{align}
  Note that for~$k=1$, \eqref{eq:matirxm_k} reduces to~\eqref{eq:mmeas}.

\begin{Remark}\label{rem:disi}
 For an integer~$k\in\{1,\dots,n\}$, the system~\eqref{eq:nonlin}
is said to be $k$-order contracting
if there exists a matrix measure~$\mu$ such that
\be\label{eq:disi}
\mu(J_f^{[k]}(t,x))\leq-\eta<0,
\ee
for any~$t\geq 0$ and any~$x$ in the state-space~\cite{kordercont}
(see also~\cite{leonov,weak_slotine}). Roughly speaking, this implies that the dynamics contracts $k$-dimensional volumes. To explain this, 
consider the system~$\dot x(t)=A(t)x(t)$, and
suppose that it  is~$n$-order contracting. Since~$A^{[n]}(t)=\trace(A(t))$, this implies that~$\trace(A(t))\leq-\eta<0$ for all~$t\geq 0$.
Combining this with~\eqref{eq:xk} implies the following. Fix~$n$ initial conditions~$a^i$, $i=1,\dots,n$, 
 and let
\begin{align*}
    X(t):= \begin{bmatrix} 
    x(t,a^1)&x(t,a^2)&\dots & x(t,a^n)\end{bmatrix}.
\end{align*}
Then $\vert \det (X(t) ) \vert \leq \exp(-\eta t)  \vert \det( X(0) )\vert$. Therefore, $n$-dimensional volumes contract
at an exponential rate. 
Following  the terminology  used in physics, we say that~\eqref{eq:nonlin}
is \emph{dissipative}   if it is $n$-order contracting.
\end{Remark}
  
  From~\eqref{eq:mmeas} and~\eqref{eq:matirxm_k}, it is straightforward to obtain the next result.
\begin{Proposition}[{\cite[Corollary~1]{kordercont}}] \label{prop:monomeak}
 Let $A \in \C^{n \times n}$ and $p \in \{1, 2, \infty\}$. Suppose that there exists an integer~$\ell \in \{1, \dots, n\}$ such that 
 \be \label{eq:mulem}
 \mu_p(A^{[\ell]}) \leq 0,
 \ee
 then
 \be \label{eq:monomeak1}
 \mu_p(A^{[\ell]}) \geq \mu_p(A^{[\ell+1]}) \geq \cdots \geq \mu_p(A^{[n]}).
 \ee
\end{Proposition}
In other words,  if~$\dot X= AX$ is~$\ell$-order contracting w.r.t.~$|\cdot|_p$ then it is also~$(\ell+1)$-order contracting,
$(\ell+2)$-order contracting, and so on~\cite{kordercont}.

The next two sections  describe our main results. 

\section{$\alpha$ compounds}\label{sec:main}
In this section, we 
define the new notions of the~$\alpha$~multiplicative and $\alpha$~additive compound of a square matrix, and analyze the properties of these compounds.
In the next section, we show how this leads to the new notion of~$\alpha$~contracting systems, with~$\alpha$ \emph{real}. 

\subsection{$\alpha$~multiplicative compound }
 Consider a non-integer real number~$\alpha \in (1, n)\setminus \mathbb{Z}$. From here on we decompose~$\alpha$  as
\[
\alpha= k+s, ~ k \in \{1, 2, \dots, n-1 \}, ~ s \in (0,1).
\]
\begin{Definition} \label{def:genm}
Let~$A \in \C^{n\times n }$ be non-singular. The~$\alpha$~multiplicative compound matrix of~$A$
is  
\begin{equation} \label{eq:defgenm}
A^{(\alpha) } :=( A ^{(k)})^{1-s} \otimes (A^{(k+1)})^s .
\end{equation}
\end{Definition}

Note that $A^{(\alpha)} \in \C^{r \times r}$, where~$r : = \binom{n}{k}\binom{n}{k+1}$, and that~$A^{(\alpha)}$ may be complex (non-real) even if~$A$ is real. Since $A$  is non-singular,    $A^{(\ell)}$ is   non-singular for all~$\ell \in \{1, \dots, n\}$, so~$(A^{(k)})^{1-s}$ and~$(A^{(k+1)})^s$ in~\eqref{eq:defgenm} are well-defined.

For example, for~$\alpha=2.5$, we have~$k=2$ and~$s=1/2$, so
$
A^{(2.5) }  = (A ^{(2)})^{1/2} \otimes (A^{(3)})^{1/2},
$
which can be interpreted as a
``multiplicative interpolation'', with equal
weights,  between~$A^{(2)}$
and~$A^{(3)}$.
\begin{Example}\label{exa:diagexa}
Suppose that~$D=\diag(d_1,\dots,d_4)$ is  non-singular. Fix~$\alpha  \in (2,3)$, so that~$k=2$ and~$s=\alpha-2 \in (0,1)$. 
Then
\begin{align*}
    D^{(\alpha)}&=( D ^{(2)})^{1-s} \otimes (D^{(3)})^s\\
    &=
   \diag( (d_1 d_2)^{1-s},(d_1 d_3)^{1-s},\dots,
    (d_3 d_4)^{1-s} )
   )\\&
   \otimes
   \diag(
      (d_1 d_2 d_3)^{s},
      (d_1 d_2 d_4)^{s}, 
      (d_1 d_3 d_4)^{s}, 
      (d_2 d_3 d_4)^{s}
   )\\
   &=
   \diag( d_1 d_2 d_3^s, d_1 d_2 d_4^s, \dots,    d_2^s d_3 d_4 ) , 
\end{align*}
so, any eigenvalue of~$D^{(\alpha)}$ is a  ``multiplicative interpolation'' between   eigenvalues of~$ D ^{(2)}$ and~$ D ^{(3)}$.
\end{Example}

\begin{Remark}
If $\alpha $ is allowed to be an integer, say, $\alpha=k$ then~$s=0$ and~\eqref{eq:defgenm}  becomes 
\[
A^{(\alpha) } = (A ^{(k)})^{1} \otimes (A^{(k+1)})^0=
A^{(k)} \otimes I_r,
\]
where~$r:=\binom{n}{k+1} $. This is not equal to~$A^{(k)}$ (but, ignoring multiplicity, 
it has  the same eigenvalues as~$A^{(k)}$). Therefore, Definition~\ref{def:genm} only considers the case where $
\alpha$ is not an integer. For the integer case, we will just use the standard definition for the $k$ multiplicative compounds.
\end{Remark}

 \begin{Remark}\label{rem:symmetric}
 Suppose that~$A\in  \C^{n\times n} $ is non-singular. 
 Using Property~\eqref{pp:realT} in Prop.~\ref{prop:real}, Property~\eqref{pp:tranp} in Lemma~\ref{lem:pp}, and the fact that $(X^{(\ell)})^T = (X^T)^{(\ell)}$ for any $X \in \C^{n \times n}$ and~$\ell \in \{1, \dots, n\}$  yields
 \begin{align*}
(A^{(\alpha) })^T  &=\left (( A ^{(k)})^{1-s} \otimes (A^{(k+1)})^s \right )^T\\
&=  (( A ^{(k)})^{1-s} )^T \otimes( (A^{(k+1)})^s ) ^T\\
&=( (A^T) ^{(k)} )^{1-s} \otimes( (A^T) ^{(k+1)} )^{s}\\
&=(A^T)^{(\alpha)},
\end{align*}
In particular, if~$A=A^T$ then~$ A^{(\alpha) } =(A^{(\alpha) })^T$.
\end{Remark}
 
An alternative possible definition of the~$\alpha$~multiplicative compound matrix is
\begin{equation} \label{eq:defgenmal}
A^{(\alpha)}_\text{alt} := (A ^{1-s})^{(k)} \otimes (A^s)^{(k+1)}.
\end{equation}
The next result shows that \eqref{eq:defgenm} and \eqref{eq:defgenmal} are equivalent.
This is useful because as we will see below 
some results
are easier to derive   using the definition in~\eqref{eq:defgenm}
and others using the alternative definition~\eqref{eq:defgenmal}.

\begin{Theorem} \label{thm:equiv}
Consider a non-singular matrix $A \in \C^{n \times n}$ and fix $\alpha \in (1, n)\setminus \mathbb{Z}$.
Then 
\begin{equation}  \label{eq:equiv}
A^{(\alpha)} =  A^{(\alpha)}_\text{alt}.
\end{equation}
\end{Theorem}

\begin{IEEEproof}
Fix~$k \in \{1,\dots ,n\}$ and  $ s \in (0,1)$. It is enough to show that
 \be\label{eq:enough}
 (A^{(k)})^s =  (A^s)^{(k)}. 
 \ee
We first
consider the case when $A$ is diagonalizable,
that is, there exist  a  non-singular~$T \in \C^{n \times n}$, and a diagonal matrix~$D \in \C^{n \times n}$ such that 
$
A = T^{-1} D T,
$
which is also the Jordan canonical form of~$A$. 
Then
$
(A^{(k)})^s = ((T^{-1})^{(k)} D^{(k)} T^{(k)})^s.
$
Using the fact that~$(T^{-1})^{(k)} = (T^{(k)})^{-1}$ and Property~\eqref{pp:realinv} in Prop.~\ref{prop:real} gives 
\[
(A^{(k)})^s = (T^{(k)})^{-1} ( D^{(k)})^s T^{(k)}.
\]
Since $D^{(k)}$ is also diagonal, $( D^{(k)})^s =( D^s)^{(k)} $. 
Hence, 
\begin{align*}
    (A^{(k)})^s &= (T^{-1})^{(k)} (D^s)^{(k)} T^{(k)}\\& = (T^{-1} D^s T)^{(k)} \\&=  (A^s)^{(k)} .
\end{align*}
We conclude that~\eqref{eq:enough}
holds when~$A$ is diagonalizable. 
The proof in the general case follows from the fact 
 that diagonalizable matrices are dense in $\C^{n \times n}$  (see, e.g., \cite[Corollary~7.3.3]{lewis1991matrix}).
\end{IEEEproof}

The next example demonstrates
Theorem~\ref{thm:equiv}.

\begin{Example}
Consider  
$
A =
\begin{bmatrix}
a & 1 & 0 \\
0 & a & 0 \\
0 & 0  & b
\end{bmatrix},
$
with~$a,b\not =0$. Note that $A$ is not diagonalizable.
A straightforward computation yields
$
A^{(2)} = 
\begin{bmatrix}
a^2 & 0 & 0 \\
0 & ab & b \\
0 & 0  & ab
\end{bmatrix}
$, and the Jordan decomposition of this matrix is
$ A^{(2)} = 
\begin{bmatrix}
1 & 0 & 0 \\
0 & b & 0 \\
0 & 0  & 1
\end{bmatrix}
\begin{bmatrix}
a^2 & 0 & 0 \\
0 & ab & 1 \\
0 & 0  & ab
\end{bmatrix}
\begin{bmatrix}
1 & 0 & 0 \\
0 & 1/b & 0 \\
0 & 0  & 1
\end{bmatrix}.
$
 Definition~\ref{def:rpsm} gives
\begin{align*}
(A^{(2)})^{ s } =& 
\begin{bmatrix}
1 & 0 & 0 \\
0 & b & 0 \\
0 & 0  & 1
\end{bmatrix}
\begin{bmatrix}
a^{2s} & 0 & 0 \\
0 & a^s b^s &  s(ab)^{s-1} \\
0 & 0  & a^s b^s
\end{bmatrix}
\begin{bmatrix}
1 & 0 & 0 \\
0 & 1/b & 0 \\
0 & 0  & 1
\end{bmatrix} \\
=&
\begin{bmatrix}
a^{2s} & 0 & 0 \\
0 & a^s b^s &  s a^{s-1}b^s \\
0 & 0  & a^s b^s
\end{bmatrix}.
\end{align*}
On the other-hand,
\begin{align*}
(A^s)^{(2)}& =
\begin{bmatrix}
a^s & sa^{s-1} & 0 \\
0 & a^s & 0 \\
0 & 0  & b^s
\end{bmatrix}^{(2)}
= 
\begin{bmatrix}
a^{2s} & 0 & 0 \\
0 & a^s b^s &  s a^{s-1}b^s \\
0 & 0  & a^s b^s
\end{bmatrix},
\end{align*}
so $(A^{(2)})^{ s } = (A^s)^{(2)}$. Since $A^{(1)} =  A$ and $A^{(3)} = \det(A)$,  this implies that $A^{(\alpha)} =  A^{(\alpha)}_\text{alt}$ for any~$\alpha \in (1, 3) \setminus \mathbb{Z}$.
\end{Example}

The following discussion  shows that, unlike the standard multiplicative compound matrix, the formula~$(AB)^{(\alpha)} = A^{(\alpha)}B^{(\alpha)}$  does not hold in general.
\begin{Remark}\label{rem:comm}
Eq.~\eqref{eq:defgenm} yields
\begin{equation}\begin{split} \label{eq:abalp} 
   (AB)^{(\alpha)} =& ((AB)^{(k)})^{1-s} \otimes ((AB)^{(k+1)})^s \\
   =& (A^{(k)} B^{(k)})^{1-s} \otimes (A^{(k+1)}B^{(k+1)})^s.
\end{split}\end{equation}
On the other-hand, using Property~\eqref{pp:mixpro} in Lemma~\ref{lem:pp} gives
\begin{equation}\begin{split} \label{eq:abalp1} 
   A^{(\alpha)}B^{(\alpha)} = & ((A^{(k)})^{1-s} \otimes (A^{(k+1)})^s ) ((B^{(k)})^{1-s} \otimes (B^{(k+1)})^s ) \\
   =& ( (A^{(k)})^{1-s} (B^{(k)})^{1-s} ) \otimes ( (A^{(k+1)})^s (B^{(k+1)})^s ) .
\end{split}\end{equation}
This shows that $(AB)^{(\alpha)} \neq A^{(\alpha)}B^{(\alpha)}$ in general.  
However it may hold in some special cases. 
Suppose for example that  both $A$ and $B$ are diagonalizable, and that~$A$ commutes with~$B$. Then, the same holds for~$A^{(\ell)},B^{(\ell)}$, $\ell\in\{1,\dots,n\}$. This implies that~$A^{(\ell)}$ and~$B^{(\ell)}$ are simultaneously diagonalizable. Then it is easy to show that for any~$s\in  (0,1)$, we have~$(A^{(\ell)} B^{(\ell)})^s=(A^{(\ell)})^s (B^{(\ell)})^s$, and this implies that in this special case~$
(AB)^{(\alpha)}=A^{(\alpha)}
B^{(\alpha)}.
$
\end{Remark}

\subsection{Spectral properties of~$A^{(\alpha)}$}
Let $\lambda_i(A) \in \C, \sigma_i(A) \in \R_+$, $i = 1, \dots, n$, denote the eigenvalues and   singular values of $A \in \R^{n \times n}$, respectively, ordered such that
\begin{align} \label{eq:ord_eig}
\text{Re}(\lambda_1(A)) \geq \text{Re}(\lambda_2(A)) \geq \cdots \geq \text{Re}(\lambda_n(A)), 
\end{align}
and
\begin{align} \label{eq:ord_sig}
\sigma_1(A) \geq \sigma_2(A) \geq \cdots \geq \sigma_n(A) \geq 0.
\end{align}

\begin{Lemma}\label{lem:eigalp}
Fix a non-singular matrix $A \in \R^{n \times n}$  and   $\alpha \in (1, n) \setminus \mathbb{Z}$. Write~$\alpha=k+s$,
with~$k$ an integer and~$s\in(0,1)$.
 \begin{enumerate}[(i)]
     \item \label{p:1}
     The eigenvalues of $A^{(\alpha)}$ are
\[
\prod_{i \in Q^{k,n}_\ell} ( \lambda_i(A))^{1-s}   \prod_{i \in Q^{k+1,n}_j} ( \lambda_i(A))^{ s}
 \]
for~$\ell\in\{1,\dots,\binom{n}{k}\}$,
 $j\in\{1,\dots,\binom{n}{k+1}\}$.
     \item \label{p:2} 
     The eigenvalues of  $(A^T A)^{(\alpha)}$ are
  \[
\prod_{i \in Q^{k,n}_\ell} ( \sigma_i(A))^{2(1-s)} 
\prod_{i \in Q^{k+1,n}_j} ( \sigma_i(A))^{2 s}
 \]
for~$\ell\in\{1,\dots,\binom{n}{k}\}$,
 $j\in\{1,\dots,\binom{n}{k+1}\}$.
     \end{enumerate}
 \end{Lemma}
 
 \begin{IEEEproof}
It is well-known~\cite{muldo1990}
that~$\eta\in\C$ is an eigenvalue of~$A^{(k)}$ iff it
is the product of~$k$
eigenvalues of~$A$, that is, there exists~$\ell \in \{1,\dots,\binom{n}{k}\}$ such that
$
        \eta=\prod_{i \in Q^{k,n}_\ell} \lambda_i(A).
$
 By Definition~\ref{def:rpsm},  
 \[
 \eta^{1-s}= \prod_{i \in Q^{k,n}_\ell} ( \lambda_i(A))^{1-s}
  \]
   is an eigenvalue of $(A^{(k)})^{1-s}$.
 Similarly, every eigenvalue of~$(A^{(k+1})^s$ has the form
 \[
  \prod_{i \in Q^{k+1,n}_j} ( \lambda_i(A))^{ s}
 \]
 for some~$j\in\{1,\dots,\binom{n}{k+1}\}$.
 Using Property~\eqref{pp:eigotimes} in Lemma~\ref{lem:pp} proves Property~\eqref{p:1}. 

 To prove  Property~\eqref{p:2}, note that
$
(A^T A)^{(\alpha)} = ((A^T A)^{(k)})^{1-s} \otimes ((A^T A)^{(k+1)})^s.
$
The eigenvalues of~$A^T A$ are
 $\sigma^2_i(A)$, $i = 1, \dots, n$, and using
  Property~\eqref{p:1} yields Property~\eqref{p:2}.
  \end{IEEEproof}

\begin{Remark}\label{rem:notl2}
Fix a real~$\alpha \geq 1$,
and define~$\omega_{\alpha}  : \R^{n\times n} \to \R_+$ by
\[
\omega_{\alpha} (A): =   \sigma_1(A)\cdots \sigma_k(A )
    \left (   \sigma_{k+1}(A) \right )^{ s} .
\] 
This function  plays a  crucial  role in the  Douady
 and Oesterl\'{e} Theorem~\cite{Douady1980}, see Section~\ref{sec:app} below.
Combining  Property~\eqref{p:2} 
with the 
  ordering of    the eigenvalues and singular values implies  that 
 \begin{align}\label{Eq:qufa}
    \lambda_1((A^T A)^{(\alpha)})  = & 
   \left  ( \sigma_1(A)\cdots \sigma_k(A) \right )^{2(1-s)} \nonumber \\  
    & \left ( \sigma_1(A)\cdots \sigma_{k+1}(A) 
    \right )^{2s}  \nonumber \\
      =&  (\omega_\alpha(A)) ^2.  
 \end{align}
Since~$A^TA$ is symmetric, it follows from Remark~\ref{rem:symmetric} that so
is~$(A^TA)^{(\alpha)}$.
Hence,~$  \lambda_1((A^T A)^{(\alpha)})  =\sigma_1((A^T A)^{(\alpha)})= || (A^T A)^{(\alpha)} ||_2 $, so we conclude that
 \begin{align} \label{eq:l2rep}
  || (A^T A)^{(\alpha)} ||_2 
      &= (\omega_\alpha(A)) ^2. 
 \end{align}
Thus, the $\alpha$ multiplicative
compound provides a matrix norm expression for~$\omega_\alpha  $.
This was our original  motivation for introducing the~$\alpha$
multiplicative compound.

Note also that in general
 \begin{align*}
(\sigma_1(A^{(\alpha)}))^2 &= \lambda_1 (  (A^{(\alpha)})^T A^{(\alpha)}   )\\
              &\not = \lambda_1 (  (A^TA ) ^{(\alpha)}).
\end{align*}
  \end{Remark}

  We now turn to defining a generalization of the $k$ additive compound. 
 \subsection{$\alpha$~additive compound  }
 The definition of the~$\alpha$~additive compound matrix
 follows~\eqref{eq:defadd}.

\begin{Definition} \label{def:gena}
Let $A \in \R^{n \times n}$ and   $\alpha \in (1, n)\setminus \mathbb{Z}$. The~$\alpha$~additive compound matrix of~$A $ is 
 \be \label{eq:defgena}
A^{[\alpha]} := \frac{d}{d\varepsilon}
 (I+\varepsilon A)^{(\alpha)} |_{\varepsilon=0}.
 \ee
\end{Definition}
 
 Note that for any $\varepsilon >0$ sufficiently small and any $k \in \{ 1, 2,\dots, n\}$,  $(I+\varepsilon A)^{(k)}$ is non-singular
 and  $(I + \varepsilon A)^{(k)} \in \Omega_n$. Hence, Proposition~\ref{prop:sufreal} and Definition~\ref{def:genm} guarantee that $A^{[\alpha]}$ is well-defined and is a  real matrix. Note also that~\eqref{eq:defgena}
 implies that 
\be\label{eq:tay}
(I+\varepsilon A)^{(\alpha)} =I+\varepsilon  A^{[\alpha]} +o(\varepsilon).
\ee

\begin{Example}
Suppose that~$D=\diag(d_1,\dots,d_4)$ is   non-singular. Fix~$\alpha  \in (2,3)$, so that~$k=2$ and~$s=\alpha-2\in(0,1)$. Let~$h_i(\varepsilon):=1+\varepsilon d_i$. 
 By Example~\ref{exa:diagexa},
\begin{align*}
(I+\varepsilon     D)^{(\alpha)}
   = &
   \diag(  h_1(\varepsilon)  h_2(\varepsilon)  h_3^s(\varepsilon) ,  h_1(\varepsilon)
    h_2(\varepsilon)  h_4^s(\varepsilon) , \\
    & \dots,    h_2^s(\varepsilon)   h_3(\varepsilon)  h_4(\varepsilon))\\
     =& I+ \varepsilon \diag(d_1+d_2+s d_3, d_1+d_2+s d_4,\\
     &\dots, s d_2+d_3+d_4
    ) +o(\varepsilon),
\end{align*}
 and~\eqref{eq:defgena} yields
 \[
     D^{[\alpha]}=\diag(d_1+d_2+s d_3, d_1+d_2+s d_4,\dots, s d_2+d_3+d_4
    ).
 \]
 Note that every eigenvalue of~$D^{[\alpha]}$
 is an ``additive interpolation'' of three eigenvalues of~$D$.
\end{Example}

The next result provides an
 expression for~$A^{[\alpha]}$ in terms of~$A^{[k]}$ and~$A^{[k+1]}$.

\begin{Theorem} \label{thm:dotgenm}
Fix $A \in \R^{n \times n}$ and   $\alpha \in (1, n)\setminus \mathbb{Z}$. Then
\be \label{eq:defgena1}
 A^{[\alpha]} = ((1-s)A^{[k]}) \oplus (sA^{[k+1]}).
\ee
\end{Theorem}

Note that this also shows that~$A^{[\alpha]}$ is   real, as $A^{[\ell]}$ is real for any~$\ell \in \{1, \dots, n\}$.

\begin{IEEEproof}
Consider the case where~$A$ is diagonalizable, that is, $
A = T^{-1} D T,
$
where $T \in \C^{n \times n}$ is non-singular, and  $D \in \C^{n \times n}$ is a diagonal matrix. Let~$B(\varepsilon):=I+\varepsilon A$.
Fix $k \in \{1, \dots, n\}$.
Then 
\begin{align*}
    (B(\varepsilon ))^{(k)}& = (T^{-1}(I + \varepsilon D) T)^{(k)} \\&=  (T^{-1})^{(k)} (I + \varepsilon D)^{(k)} T^{(k)}.
\end{align*}
Since~$D$ is diagonal, so is
 $(I + \varepsilon D)^{(k)}$. Therefore, $(B( \varepsilon  ))^{(k)} $ satisfies the conditions in Lemma~\ref{lem:dotAalp}. 
 We use Lemma~\ref{lem:dotAalp} to determine the derivative of~$(B(\varepsilon))^{(\alpha)}$ with respect to~$\varepsilon$. To simplify the notation, we write~$B$ for~$B(\varepsilon)$. Then
\begin{align*}
 \frac{d}{d\varepsilon}B^{(\alpha)} =& 
 \frac{d}{d\varepsilon} \left ( (B^{(k)})^{1-s} \otimes (B^{(k+1)})^s\right  ) \\
 =&  \left( \frac{d}{d\varepsilon}  (B ^{(k)})^{1-s} \right) \otimes (B^{(k+1)})^s  \\
 &+  (B^{(k)})^{1-s} \otimes  \left( \frac{d}{d\varepsilon} (B^{(k+1)})^s \right)\\
 =& \left((1-s) (B^{(k)})^{-s} \frac{d}{d\varepsilon}  B^{(k)} \right) \otimes (B^{(k+1)})^s \\
 &+ (B^{(k)})^{1-s} \otimes \left(s (B^{(k+1)})^{s-1} \frac{d}{d\varepsilon} B^{(k+1)} \right).
\end{align*}
Setting $\varepsilon = 0$ and using the fact that~$B(\varepsilon)|_{\varepsilon=0}=I$ and~\eqref{eq:defadd}  yields
\begin{align*}
  \frac{d}{d\varepsilon}(I+\varepsilon A)^{(\alpha)} |_{\varepsilon=0} =& 
 (1-s) A^{[k]} \otimes I_{r_1} + I_{r_2} \otimes (s  A^{[k+1]} ) \\
 =& ((1-s)A^{[k]}) \oplus (sA^{[k+1]}),
\end{align*}
where $r_1 : =\binom{n}{k+1}$, and $r_2 : =\binom{n}{k}$. 
This completes the proof when~$A$ is  diagonalizable, and the general case follows from using a similar argument as in the proof of Theorem~\ref{thm:equiv}. 
\end{IEEEproof}

 \begin{Remark}\label{rem:symmetric_new}
 Suppose that~$A\in  \R^{n\times n} $. 
 Then it is easy to see that Theorem~\ref{thm:dotgenm}
 implies that
$
 (A^{[\alpha]})^T=(A^T)^{[\alpha]}.
 $ 
 In particular, if~$A$ is symmetric then so is~$A^{[\alpha]}$. 
\end{Remark}

The next result shows that the $\alpha$ additive compound  satisfies the  same additivity
property as the $k$ additive compound.

\begin{Theorem} \label{thm:addp}
Let $A , B \in \R^{n \times n}$ and   $\alpha \in (1, n)\setminus \mathbb{Z}$. Then, 
\[
( A + B)^{[\alpha]} = A^{[\alpha]} + B^{[\alpha]}.
\]
\end{Theorem}

\begin{IEEEproof}
Using~\eqref{eq:defgena1}
gives 
\begin{align*}
   ( A + B)^{[\alpha]} = &((1-s)(A + B)^{[k]}) \oplus (s(A + B)^{[k+1]}) \\
   =& (1-s)(A + B)^{[k]} \otimes I_{r_1} +   s I_{r_2} \otimes (A + B)^{[k+1]}.
\end{align*}
Applying~\eqref{eq:akadd} and  Properties~\eqref{pp:asos}-\eqref{pp:dist1} in Lemma~\ref{lem:pp} yields
\begin{align*}
   ( A + B)^{[\alpha]}  
   =& (1-s)A ^{[k]} \otimes I_{r_1} +   s I_{r_2} \otimes A^{[k+1]}  \\
   &+ (1-s)B ^{[k]} \otimes I_{r_1} + s I_{r_2} \otimes B^{[k+1]} \\
   =& ((1-s)A^{[k]}) \oplus (sA^{[k+1]}) \\& + ((1-s)B^{[k]}) \oplus (sB^{[k+1]}) \\
   =& A^{[\alpha]} + B^{[\alpha]},
\end{align*} 
and this completes the proof.
\end{IEEEproof}

The next result shows that~\eqref{eq:expat} also holds for the $\alpha$  compounds.  

\begin{Theorem} \label{thm:expgen}
Let $A \in \R^{n \times n}$ and   $\alpha \in (1, n) \setminus \mathbb{Z}$. Then, 
\be  \label{eq:expgen}
\exp(A^{[\alpha]}t) = (\exp(A t))^{(\alpha)}, \text{ for any } t \in \R.
\ee
\end{Theorem}

\begin{IEEEproof}
Using the alternative definition~\eqref{eq:defgenmal} gives 
\begin{align*}
  (\exp(At))^{(\alpha)}_\text{alt} =& ((\exp(At))^{1-s})^{(k)}\otimes ((\exp(At))^{s})^{(k+1)} \\
  =& (\exp(At(1-s)) )^{(k)}\otimes (\exp(Ats))^{(k+1)} , 
\end{align*}
and applying~\eqref{eq:expat} and
 Properties~\eqref{pp:asos} and \eqref{pp:expkro} in Lemma~\ref{lem:pp},
 gives
 \begin{equation}\begin{aligned}
  (\exp(At))^{(\alpha)}_\text{alt}  
  =& \exp(A^{[k]}t(1-s)) \otimes \exp(A^{[k+1]}ts) \\
  =& \exp \left ((A^{[k]}t(1-s))\oplus(A^{[k+1]}ts) \right ) \\
  =& \exp\left (((A^{[k]}(1-s))\oplus(A^{[k+1]}s)) t \right ) \\
  =& \exp(A^{[\alpha]}t).
\end{aligned}\end{equation}
Applying Theorem~\ref{thm:equiv} completes the proof.
\end{IEEEproof}

Combining Lemma~\ref{lem:eigalp}, Theorem~\ref{thm:expgen}, and the fact that the eigenvalues of  $\exp(At)$ are
 $\exp(\lambda_i(A)t)$, $i\in\{1,\dots,n\}$, yields  the following result.
 
\begin{Corollary}
\label{coro:eigalpa}
Let $A \in \R^{n \times n}$ and   $\alpha \in (1, n)\setminus \mathbb{Z}$.
 The eigenvalues of $A^{[\alpha]}$ are
 \[
  (1-s)\sum_{i \in Q^{k,n}_\ell} \lambda_i(A) +  s \sum_{i \in Q^{k+1,n}_j} \lambda_i(A)
 \]
for~$\ell\in\{1,\dots,\binom{n}{k}\}$,
 $j\in\{1,\dots,\binom{n}{k+1}\}$.
 \end{Corollary}
 
 This implies in particular that~$\lambda_1(A^{[\alpha]})$, that is, the eigenvalue of~$A^{[\alpha]}$ with the largest real part is
 \begin{align} \label{eq:specalpha}
      \lambda_1(A^{[\alpha]})& = (1-s)\sum_{i=1}^k \lambda_i(A) + s\sum_{i=1}^{k+1}\lambda_{i}(A) \nonumber \\
      &= \sum_{i=1}^k \lambda_i(A) + s \lambda_{k+1}(A). 
  \end{align}

 \begin{Remark}
 By the spectral properties of the standard additive compounds \cite{muldo1990}, $\sum_{i \in Q^{k,n}_\ell} \lambda_i(A)$ and $\sum_{i \in Q^{k+1,n}_j} \lambda_i(A)$ are the eigenvalues of $A^{[k]}$ and $A^{[k+1]}$, respectively. Now    Property~\eqref{pp:eigoplus} in Lemma~\ref{lem:pp}  implies that~$(1-s)\sum_{i \in Q^{k,n}_\ell} \lambda_i(A) +  s \sum_{i \in Q^{k+1,n}_j} \lambda_i(A)$ is an eigenvalue of $A^{[\alpha]}$. This provides another proof for Corollary~\ref{coro:eigalpa}.
 \end{Remark}

 The next result     follows from Thm.~\ref{thm:expgen}.
 
 \begin{Corollary} \label{coro:tansi}
 Suppose that $X(t)$ is the solution of the linear time-invariant matrix differential equation
 \be \label{eq:xta}
 \frac{d}{dt}{X}(t) = A X(t), ~X(0) = I, 
 \ee
 with $A \in \R^{n \times n}$. Fix~$\alpha \in (1,n) \setminus \mathbb{Z}$. Then  
 \be \label{eq:xalphax}
 \frac{\diff}{\diff t} X^{(\alpha)}(t) = A^{[\alpha]} X^{(\alpha)}(t), \quad X^{(\alpha)}(0)=I.
 \ee
 \end{Corollary}
 \begin{IEEEproof}
The solution of~\eqref{eq:xta} is $X(t) = \exp(At)$. By~\eqref{eq:expgen}, $X^{(\alpha)}(t) = \exp(A^{[\alpha]}t)$, and this yields~\eqref{eq:xalphax}.
 \end{IEEEproof}
 
 It is important to note that
Eq.~\eqref{eq:xalphax}     holds  only when~$A$ is a  constant matrix.   The next example demonstrates this.
 
 \begin{Example}
 Consider the linear time-varying matrix differential equation
 \[
 \dot{X}(t) = A(t) X(t),\quad X(0) = I ,
 \]
 with 
$
 A(t) =
 \begin{bmatrix}
  1 & 0 \\
 1 & t 
 \end{bmatrix}.
$
 In this case,  
$
 X(t) =
 \begin{bmatrix}
 a(t) & 0 \\
 b(t) & c(t)
 \end{bmatrix}
$ with \begin{align*}
  a(t)& := \exp(t),~
  b(t) := \int_0^t \exp\left (\sigma + \frac{t^2 -\sigma^2}{2} \right ) \diff \sigma,\\
  c(t)&:= \exp( {t^2}/ 2).
\end{align*}    For any~$t\in(0,\infty)\setminus\{2\}$, 
the Jordan decomposition of~$X$ is 
\[
X  = 
\begin{bmatrix}
\frac{a-c}{b} & 0 \\
1 & 1 
\end{bmatrix}
\begin{bmatrix}
a & 0 \\
0 & c 
\end{bmatrix}
\begin{bmatrix}
\frac{b}{a-c} & 0 \\
-\frac{b}{a-c} & 1 
\end{bmatrix},
\]
so 
\[
X^{s} = 
\begin{bmatrix}
\frac{a-c}{b} & 0 \\
1 & 1 
\end{bmatrix}
\begin{bmatrix}
a^s & 0 \\
0 & c^s
\end{bmatrix}
\begin{bmatrix}
\frac{b}{a-c} & 0 \\
-\frac{b}{a-c} & 1 
\end{bmatrix}
=
\begin{bmatrix}
a^s & 0 \\
-\frac{a^s -c^s}{a-c}b & c^s 
\end{bmatrix}.
\]
Let $\alpha = 1+s$, with~$s\in(0,1)$. Then Definition~\ref{def:genm}
gives
$ X^{(\alpha)} = \left(ac \right)^s X^{1-s}$,
and~\eqref{eq:defgena1} gives
$ A^{[\alpha]} =  
\begin{bmatrix}
1 + st & 0 \\
1-s & s+t
\end{bmatrix}
$.
In particular, substituting~$s=0.5$ and $t = 1$ yields 
\begin{align*}
 \frac{\diff }{\diff t} X^{(1.5)}(1) = & 
\begin{bmatrix}
5.23551 & 0 \\
4.48053 & 4.07742
\end{bmatrix} \\
A^{[1.5]}(1) X^{(1.5)}(1) =& 
\begin{bmatrix}
5.23551 & 0 \\
4.26352 & 4.07742
\end{bmatrix} ,
\end{align*}
so $\frac{\diff }{\diff t} X^{(\alpha)}(t) \not = A^{[\alpha]}(t) X^{(\alpha)}(t)$.
\end{Example}

 It is well-known~\cite{muldo1990}
 that for any~$A,T \in \C^{n\times n}$, with~$T$ non-singular, and  any integer~$\ell\in\{1,\dots,n\}$, we have
\be\label{eq:linttra}
(TAT^{-1})^{[\ell]}=T^{(\ell)} 
A^{[\ell]} (T^{-1})^{(\ell)} = T^{(\ell)} 
A^{[\ell]} (T^{(\ell)})^{-1}.
\ee
The next result shows 
the $\alpha$~additive compounds
under a coordinate transformation.
\begin{Theorem} \label{thm:coortans}
Let~$A,T \in \C^{n\times n},
$ with~$T$ non-singular, and 
pick~$\alpha\in(1,n)\setminus \mathbb{Z}$.
Then
\be \label{eq:coortran}
(TAT^{-1})^{[\alpha]}=(T^{(k)} \otimes T^{(k+1)}) A^{[\alpha]}(T^{(k)} \otimes T^{(k+1)})^{-1}.
\ee
\end{Theorem}
\begin{IEEEproof}
Let~$B:=TAT^{-1}$. 
Using \eqref{eq:defgena1} and~\eqref{eq:linttra} yields
\begin{align*}
 B^{[\alpha]}  
 =& ((1-s)(TAT^{-1})^{[k]}) \oplus (s(TAT^{-1})^{[k+1]})\\
=&
 (1-s)( T^{(k)} A^{[k]}  (T^{-1})^{(k)}) 
 \otimes (T^{(k+1)}   I_{r_1} (T^{(k+1)})^{-1} ) \\ 
 & +
 s  (T^{(k)} I_{r_2} (T^{(k)})^{-1}) \otimes
 (  T^{(k+1)} A^{[k+1]}  (T^{-1})^{(k+1)}).
\end{align*}
Using Properties~\eqref{pp:mixpro} and~\eqref{pp:inverse} in Lemma~\ref{lem:pp} gives
\begin{align*}
B^{[\alpha]}  
=& (1-s)(T^{(k)} \otimes T^{(k+1)})  (A^{[k]} \otimes I_{r_1}) \\
& \left((T^{(k)})^{-1} \otimes (T^{(k+1)})^{-1} \right) 
+ s(T^{(k)} \otimes T^{(k+1)}) \\& (I_{r_2} \otimes A^{[k+1]}) 
 \left((T^{(k)})^{-1} \otimes (T^{(k+1)})^{-1} \right) \\
=& (T^{(k)} \otimes T^{(k+1)})   \left( (1-s) (A^{[k]} \otimes I_{r_1}) + sI_{r_2} \otimes A^{[k+1]} \right) \\ & (T^{(k)} \otimes T^{(k+1)})^{-1} \\
=& (T^{(k)} \otimes T^{(k+1)}) A^{[\alpha]}(T^{(k)} \otimes T^{(k+1)})^{-1},
\end{align*}
and this completes the proof.
\end{IEEEproof}

\begin{Remark}
For a non-singular $T \in \R^{n \times n}$, note that 
$
T^{(k+\frac12)} = (T^{(k)})^\frac12 \otimes (T^{(k+1)})^\frac12
$,
and thus
\begin{align*}
(T^{(k+\frac12)})^2 =& \left((T^{(k)})^\frac12 \otimes (T^{(k+1)})^\frac12 \right)\left((T^{(k)})^\frac12 \otimes (T^{(k+1)})^\frac12 \right) \\
=&\left((T^{(k)})^\frac12 (T^{(k)})^\frac12 \right) \otimes \left((T^{(k+1)})^\frac12  (T^{(k+1)})^\frac12 \right) \\
=& T^{(k)} \otimes T^{(k+1)}.
\end{align*}
Therefore, \eqref{eq:coortran} can be rewritten in  the more compact form
\[ \label{eq:coortran1}
(TAT^{-1})^{[\alpha]}=(T^{(k+\frac12)})^2 A^{[\alpha]}(T^{(k+\frac12)})^{-2}.
\]
\end{Remark}

The next subsection analyzes  the 
matrix measure 
of the $\alpha$~additive compound.  This will play an important role in the analysis of~$\alpha$ contracting systems.
\subsection{Matrix measures of the $\alpha$~additive compound}

It is well-known~\cite{mono_norms} 
that if~$|\cdot|$ is monotonic then 
the induced matrix norm satisfies
$
||D||=\max \big(|d_{1}|,\dots,|d_{n}| \big),
$
for any diagonal matrix~$D=\diag(d_1,\dots,d_n)\in\C^{n\times n}$.   This implies that
  the induced matrix measure satisfies
\begin{align}\label{eq:diagD}
  \mu(D)&=\lim_{\varepsilon \downarrow 0} \varepsilon^{-1} (||I+\varepsilon D||-1)\nonumber\\
  &=\lim_{\varepsilon \downarrow 0} \varepsilon^{-1} (\max_i\{|1+\varepsilon d_i|\}-1)\nonumber\\
  &=\max_i \{\real (d_i)\}.
\end{align}
Therefore, for any $\ell \in \{1, \dots, n\}$, we have
\begin{align}\label{eq:diagDell}
  \mu(D^{[\ell]})
  =\max_{\{i_1, \dots,i_\ell \} \in Q^{\ell, n}} \left ( \sum_{p=1}^\ell \real (d_{i_p})\right) .
\end{align}
Eqs.~\eqref{eq:defgena1} and \eqref{eq:diagDell} imply  that 
\begin{align}\label{eq:diagDalp}
  \mu(D^{[\alpha]})
  =(1-s)\mu(D^{[k]}) + s \mu(D^{[k+1]}) .
\end{align}

Our next goal is to show that this holds for general  matrices. 
Towards this end, we first provide  a useful expression for the matrix measure of a Kronecker sum of matrices.
\begin{Theorem} \label{thm:muxplusy}
Let $\mu$ denote a matrix measure associated with a induced matrix norm $\Vert \cdot \Vert$ such that 
\be \label{eq:kronorm}
\Vert A \otimes B \Vert = \Vert A \Vert \Vert B \Vert
\ee
for any matrices $A ,B $. 
Then
\be \label{eq:xyequ}
\mu(X \oplus Y) = \mu(X) + \mu(Y),
\ee
 for any $X \in \R^{n \times n}$ and~$Y\in\R^{m\times m}$.
\end{Theorem}

\begin{IEEEproof}
Fix~$\varepsilon>0$.
Properties~\eqref{pp:asos} and~\eqref{pp:expkro} in Lemma~\ref{lem:pp} yield 
\begin{align*}
  || \exp(\varepsilon ( X\oplus Y))|| &= || \exp(\varepsilon  X     \oplus \varepsilon  Y )   ||\\
   &= || \exp(\varepsilon  X )  \otimes  \exp(\varepsilon  Y )   || .
\end{align*}
By~\eqref{eq:kronorm},~$|| \exp(\varepsilon ( X\oplus Y))||  
   = || \exp(\varepsilon  X ) || || \exp(\varepsilon  Y )   ||$,
   and thus
\begin{align*}
  \frac{d}{d\varepsilon} || \exp(\varepsilon ( X\oplus Y))||  
   &=\left( \frac{d}{d\varepsilon} || \exp(\varepsilon  X ) || \right) || \exp(\varepsilon  Y )   || \\
   &+   || \exp(\varepsilon  X ) || \left(\frac{d}{d\varepsilon}  || \exp(\varepsilon  Y   )   ||\right ). 
\end{align*}
It follows from~\eqref{eq:defmu} that  
\[
\mu(A)=\left (\frac{d}{d\varepsilon} || \exp(\varepsilon A)||\right) |_{\varepsilon =0},
\]
for any~$A\in\R^{n\times n}$.
Thus,
\begin{align*}
& \mu(X \oplus Y)= 
  \left(\frac{d}{d\varepsilon} || \exp(\varepsilon ( X \oplus Y))|| \right) |_{\varepsilon =0} \\
  =& \left( \frac{d}{d\varepsilon} || \exp(\varepsilon  X ) || \right)|_{\varepsilon =0}   
   +    \left(\frac{d}{d\varepsilon}  || \exp(\varepsilon  Y )     || \right) |_{\varepsilon =0}\\
   =&\mu(X)+\mu(Y),
\end{align*}
and this completes the proof.
\end{IEEEproof}

We can now provide a useful expression for the matrix measure of~$A^{[\alpha]}$. 
\begin{Corollary}\label{coro:ineq}
Let $\mu_p$ denote a matrix measure induced by some~$L_p$ norm with~$p \geq 1$. For any~$A \in \R^{n \times n}$ and any $\alpha \in (1, n)\setminus \mathbb{Z}$, 
we have 
\be \label{eq:mupequa}
\mu_p (A^{[\alpha]}) = (1-s) \mu_p (A^{[k]}) + s\mu_p (A^{[k+1]}).
\ee
\end{Corollary}

\begin{IEEEproof}
From Prop.~\ref{prop:normkrop}, \eqref{eq:kronorm} holds for all $L_p$ norms.
Since~$ A^{[\alpha]} 
=   ( (1-s)A^{[k] }) \oplus ( s A^{[k+1]}) $,  
Theorem~\ref{thm:muxplusy} yields
\begin{align*}
\mu_p(A^{[\alpha]}) 
 =& \mu_p \left ( (1-s)A^{[k]} \right ) + \mu_p \left (  sA^{[k+1]} \right ) \\
=& (1-s)\mu_p \left( A^{[k]} \right ) + s\mu_p \left (  A^{[k+1]} \right ),
\end{align*}
where the last equality follows from the homogeneity of the matrix measure, and the fact that~$s\in(0,1)$.
\end{IEEEproof}

 The next example demonstrates Corollary~\ref{coro:ineq} in the case~$n=2$.
 \begin{Example}
 Let~$A  \in\R^{2\times 2} $.
 Fix~$\alpha\in(1,2)$. Then~$\alpha=k+s$, with~$k=1$ and~$s  \in (0,1)$, so 
\begin{align*}
    A^{[\alpha]} &= ((1-s)A^{[1]}) \oplus (sA^{[2]})\\
     &= ((1-s)A ) \oplus (s \trace(A))\\
    &= ((1-s)A)   \otimes I_1 + I_2 \otimes  (s \trace (A))\\
    &=(1-s)A+   s \trace(A) I_2  . 
\end{align*} 
Recall that for any matrix measure~$\mu$ and any~$c\in\R$, 
 $\mu(A+c I)=\mu(A)+ c$ (see e.g.~\cite{feedback2009}). Thus, 
\begin{align*}
   \mu( A^{[\alpha]} )  &=  \mu((1-s)A)+s\trace(A)\\
                        &=   (1-s)\mu(A^{[1]})+ s \mu (A^{[2]})  . 
\end{align*} 
Note that for this particular example,~\eqref{eq:mupequa} holds for all matrix measures.
 \end{Example}

 In the remainder of this paper  we always assume that~$\mu$ is induced from some~$L_p$ norm, with~$p\geq 1$.
  
 The next section
 describes an application of the~$\alpha$~compounds in the context of the 
  Douady
 and Oesterl\'{e} Theorem~\cite{Douady1980}. For a modern treatment of this theorem and its numerous extensions and applications, see the recent monograph~\cite{book_volker2021}. Some connections between contracting
 systems and the Douady and  Oesterl\'{e} Theorem 
 have already appeared in the note~\cite{weak_slotine}.

 \section{An application: $\alpha$ contracting systems}~\label{sec:app}
 
 In this section,~$\alpha\in[1,n)$,
 and the special case where~$\alpha$ is an integer is also allowed. 
 The Hausdorff dimension of a set~$K\subset\R^n$ is denoted by~$\dim_H K$. 
  Let~$D\subseteq \R^n$ be an open set, and let~$g:D \to \R^n$ be a~$C^1$ mapping.
  Let
  \[
  J_g(x):=\frac{\partial}{\partial x}g(x).
  \]
 A set~$K\subseteq D$ is said to be \emph{negatively invariant under~$g$}
 if~$K\subseteq g(K)$.  Intuitively speaking,~$g$
 ``increases''~$K$.
 The next result is the    Douady
 and Oesterl\'{e} theorem~\cite{Douady1980}. 
 We state in the form given in~\cite{smith_hauss}.
 
 \begin{Theorem}\label{thm:bwdk}
 Suppose that~$K\subset D$ is   compact  and  negatively invariant under~$g$.
 Fix~$\alpha \in [1,n)$, and write~$\alpha=k+s$,
 with~$k$ an integer and~$s\in[0,1)$.
 Let
  \be \label{eq:wds1}
  \omega (K,\alpha,g):=\max_{x\in K}\left( \sigma_1(J_g(x))\cdots \sigma_k(J_g(x))  (\sigma_{k+1}(J_g(x)))^s\right) .
  \ee
  If~$\omega (K,\alpha,g)<1$
  then
$
 \dim_H K <\alpha .
 $
 \end{Theorem}
Intuitively speaking, 
\eqref{eq:wds1} implies that 
 $g$  is  a ``contraction in dimension~$\alpha$'', uniformly in~$K$.
 If~$g$ ``increases''~$K$  
then necessarily~$\dim_H K<\alpha$.

The next simple example demonstrates Thm.~\ref{thm:bwdk}.
\begin{Example}
Let~$g:\R^3\to\R^3$ be  the linear mapping given by~$g(x)=\diag( 1,1/2,1/4  )x$. Then the singular values of the Jacobian of~$g$ are~$1,1/2,1/4$,
and~\eqref{eq:wds1} holds for any~$\alpha>1$. Thus, Thm.~\ref{thm:bwdk} implies that for any compact set~$K\subset \R^3$ such that~$K\subseteq g(K)$, we have~$\dim_H K\leq 1$. For example,  the set~$K :=[0,1]\times \{0\}\times\{0\}$ satisfies~$K\subseteq g(K)$  and~$\dim_H K = 1$. 
\end{Example}

Using the $\alpha$~multiplicative compound  we can express condition~\eqref{eq:wds1} in a more elegant form. Indeed,
it follows from~\eqref{eq:l2rep} that
\[
 ( \omega (K,\alpha,g))^2=\max_{x\in K}   
 ||\left ( J_g^T(x) J_g(x) \right )^{(\alpha)}   ||_2
  ,
\]
so the condition in Thm.~\ref{thm:bwdk} 
becomes
\[
\max_{x\in K}  
 || ( J_g^T(x) J_g(x) )^{(\alpha)}   ||_2
   <1.
\]
This provides a more intuitive description for    ``contraction in dimension~$\alpha$'' of a mapping~$g$.

Thm.~\ref{thm:bwdk}
has been used to upper bound the 
Hausdorff dimension of invariant sets (and, in particular, attractors)  of dynamical systems. Our results allow to restate and generazlize these results in a more intuitive  fashion using the~$\alpha$ additive compound.

  Consider the time-varying dynamical system:
  \be\label{eq:dsys}
  \dot x=f(t,x),
  \ee
  where~$f$ is~$C^1$. 
  Let~$x(t,t_0,x_0)$ denote the solution of~\eqref{eq:dsys} at time~$t$ with~$x(t_0)=x_0$. We assume from here on that~$t_0=0$, and let~$x(t,x_0):=x(t,0,x_0)$. We also assume that there exists an invariant set~$D\subseteq \R^n$, that is, for any~$x_0\in D$ we have~$x(t,x_0) \in D$ for  all~$t\geq 0$. Let
  $
  J_f(t,x):=\frac{\partial }{\partial x}f(t,x),
 $
  and
  consider the matrix differential equation
  \[
  \dot X(t)=J_f(x(t,x_0)) X(t), \quad X(0)=X_0.
  \]
From now on, we always consider the $L_p$ norms with $p \geq 1$ and the associated matrix measure $\mu$. We begin with an auxiliary result.
  \begin{Proposition} \label{prop:gammajf}
  Let~$K\subset\R^n$   be a  compact      invariant set of~\eqref{eq:dsys}. 
  Fix~$\alpha\in [1,n)$
  and let~$\alpha=k+s$, with~$k$   integer and~$s\in [ 0,1)$. 
 For an induced matrix measure~$\mu$ and~$t \geq 0$, let
  \begin{align*} 
\gamma_{J_f} (t) :=  &  \max_{ x_0 \in K}  \int_0^t
  \mu(J_f^{[\alpha]}( x(\tau,x_0)  )) 
  \diff \tau .
  \end{align*}
  Then
 \begin{align*}
&||X^{(k)} (t)|| ^{1-s} ||X^{(k+1)} (t ) ||^s \\
&~~\leq \exp( \gamma(t) )||X_0^{(k)}  || ^{1-s} ||X^{(k+1)}_0 ||^s, \text{ for any }x_0\in K.
\end{align*}
\end{Proposition}

  \begin{IEEEproof}
  Pick~$\ell\in\{1,\dots, n\}$.
  Since~$\frac{d}{dt} X^{(\ell)}=J^{[\ell]} X^{ (\ell )}$, 
 $
  ||X^{(\ell)} (t) || \leq \exp (\int_0^ t\mu (J^{[\ell]} (\tau)  ) \diff \tau ) ||X_0^{(\ell)}||
   .
  $
Applying this bound to~$||X^{(k)} (t)|| ^{1-s} ||X^{(k+1)} (t ) ||^s $,
and using~\eqref{eq:mupequa}   completes the proof. 
\end{IEEEproof}

  We say that a constant 
  set~$K\subseteq D$ is a strongly invariant set of~\eqref{eq:dsys}
  if
  \be\label{eq:invtau} 
  K =  x(t, K) \text{ for all } t\geq 0.
  \ee
  For example, an equilibrium or a limit cycle are strongly  invariant  sets. 
  
  We can now bound the Hausdorff dimension of strongly invariant sets of~\eqref{eq:dsys}, thus extending a result in~\cite{smith_hauss}.
 For generality, contraction theory typically uses contraction metrics~\cite{LOHMILLER1998683} and associated scaled norms. 
  Consider  a~$C^1$ scaling matrix~$\Theta: K\to \R^{n\times n}$ satisfying
  \be\label{eq:scathe}
  \det( \Theta(z)  ) \not =0   \text{ for all } z\in K.
  \ee
  Let~$ \Theta _f(z)$ denote  the matrix obtained by replacing every entry~$ \theta_{ij}(z)$ in~$ \Theta (z)$ 
by the value~$(\frac{\partial  \theta_{ij}(z)}{\partial z} )^T f(z)$, and define the so-called generalized Jacobian~\cite{LOHMILLER1998683} as
\[
\bar J:=
\Theta _f  \Theta ^{-1}  + \Theta  J_f  \Theta ^{-1}  .
\]
Note that if~$\Theta(z)=I$ for all~$z$ then~$\bar J=J_f$.
The next result bounds the Hausdorff dimension of a strongly invariant set using the generalized  Jacobian~$\bar J$. 
  \begin{Theorem}   \label{thm:rusgens}   
Let~$K\subset\R^n$   be a  compact  and strongly invariant set of~\eqref{eq:dsys}. 
  Fix~$\alpha\in [1,n)$
  and let~$\alpha=k+s$, with~$k$   integer and~$s\in [ 0,1)$. 
 Assume there exist an induced matrix measure~$\mu :\R^{n\times n}\to \R$ and~$\tau>0$ such that
  \begin{align}\label{eq:alpdef}
\gamma_{\bar J}(\tau) <0.
  \end{align}
  Then
 $
  \dim_H K<  \alpha.
$
\end{Theorem}

\begin{IEEEproof}
Define~$g:\R_+\times K\to K$ by~$g(t,x_0):=x(t,x_0)$.
Then~$J_g(t,x_0):= \frac{\partial}{\partial x_0} x(t,x_0)$. 
Let
\be \label{eq:defY}
Y(t,x_0):= \Theta  (x(t,x_0))
 \frac{\partial }{\partial x_0}g(t,x_0).
\ee
To simplify the notation, we sometimes write $\Theta(x)$ or $\Theta(t)$ for $\Theta(x(t, x_0))$.
By~\eqref{eq:defY},
\[
\dot Y \ = \ \dot \Theta  \frac{\partial }{\partial x_0}g +
  \Theta  \frac{\partial}{\partial x_0} \dot x 
\ = \ (\Theta _f J_g+ \Theta  J_f  \Theta ^{-1}) Y.
\]
 Thus,~$Y(t,x_0)$ is the solution at time~$t$ of the matrix differential equation
$
\dot Y=\bar J Y,
$
initialized with~$Y(0)= \Theta (x_0)$.
Let~$c(x_0) := ||\Theta^{(k)}(x_0) ||^{1-s} ||\Theta^{(k+1)}(x_0) ||^s $.
Prop.~\ref{prop:gammajf} and~\eqref{eq:alpdef} imply that 
$$
||Y^{(k)}(\tau)||^{1-s }||Y^{(k+1)}(\tau)||^{s }\leq  c (x_0) \exp(\gamma_{\bar J} (\tau))    ,
$$
for any~$x_0\in K$. Hence, for any integer~$\ell \geq 1$, 
\begin{align} \label{eq:thetaltau}
||\Theta^{(k)}( \ell \tau  ) J_g^{(k)}(\ell \tau  ) ||  ^{1-s}& ||\Theta^{(k+1)}( \ell \tau  ) J_g^{(k+1)}(\ell \tau  ) ||  ^{ s}  \nonumber \\& \leq  c(x_0)  \exp(\gamma_{\bar J}( \ell \tau)). 
\end{align}
 Recall that if~$||\cdot||: \C^{n \times n} \to \R_+$ is an induced matrix norm and~$P\in\C^{n\times n}$ 
 is non-singular, then the~$P$-weighted induced matrix norm is~$ ||M||_P:=|| PMP^{-1}||$. Eq.~\eqref{eq:thetaltau} yields
\begin{align*}
&||  J_g^{(k)}(\ell \tau  )    ||_{\Theta^{(k)}( \ell \tau  )}  ^{1-s} 
||J_g^{(k+1)}(\ell \tau  )  ||_{\Theta^{(k+1)}( \ell \tau  ) }  ^{ s} \\  &  \leq c(x_0)\exp(\gamma_{\bar J} (\ell \tau)) || (\Theta^{(k)}( \ell \tau  ))^{-1} ||^{1-s} ||(\Theta^{(k+1)}( \ell \tau  ))^{-1}  ||^s  .
\end{align*}
Since~$K$ is compact, 
 we can make the right-hand side of this equation arbitrarily small   
by taking~$\ell$ large enough.
Using the equivalence of norms implies that there exists an integer~$\bar \ell$ such that
$
||  J_g^{(k)}(\bar \ell \tau  )    ||_2 ^{1-s}
||J_g^{(k+1)}(\bar \ell \tau  )  ||_2 ^{ s}      <  1.
$
Let~$\sigma_i$, $i=1,\dots,n$,
denote the singular values of~$ J_g(\bar \ell \tau )$.
Then we conclude that
\[
\sigma_1\dots \sigma_k\sigma_{k+1}^s<1.
\]
Since~$g(\bar \ell \tau,K)=K$, Thm.~\ref{thm:bwdk}
implies that~$\dim_H K<\alpha$. 
\end{IEEEproof}

Of course, a sufficient condition for~\eqref{eq:alpdef} to hold is that~$
\mu(\bar J^{[\alpha]}(x))<0$    for all~$ x\in K$.

From now on we consider  for simplicity the non-scaled case, i.e.~$\bar J=J_f$. 
Thm.~\ref{thm:rusgens}   naturally leads to the following new definition. 

\begin{Definition}
Let~$\mu$ be a matrix measure induced by an~$L_p$ norm with~$p\geq 1$.
Suppose that the trajectories of~\eqref{eq:dsys} evolve on a state space~$D$. Pick a real~$\alpha\geq 1$.
System~\eqref{eq:dsys}
is called $\alpha$~contracting w.r.t. the norm~$|\cdot|_p$ if
\be \label{eq:alphacon}
\mu(J_f^{[\alpha]}(t,x))\leq -\eta<0,\text{ for all } t\geq0,\;x\in D.
\ee
\end{Definition}

 \begin{Remark}\label{rem:adding}
 An important property of contracting systems is that 
various compositions of contracting  systems yield a contracting system~\cite{LOHMILLER1998683,modular,netwok_contractive,weak_slotine}.  The subadditivity of the matrix measure and Thm.~\ref{thm:addp} suggest that this remains valid for interconnections of~$\alpha$ contracting systems. 
 As a simple example, consider the interconnected  system
 \be\label{eq:inter}
 \dot x(t)=c_1(t)f(t,x)+c_2(t)g(t,x),
 \ee
 with~$c_i(t)\geq 0$ for any~$t\geq 0$.
 The Jacobian of this system is~$c_1 J_f+c_2 J_g$,
 and
 \begin{align*}
 \mu((c_1 J_f+c_2 J_g)^{[\alpha]})&=\mu(c_1 J_f^{[\alpha]}+
 c_2 J_g ^{[\alpha]}) \\& \leq
 c_1\mu(J_f^{[\alpha]})+c_2 \mu( J_g ^{[\alpha]}).
 \end{align*}
 Thus, it is straightforward to provide sufficient conditions for~$\alpha$ contraction of~\eqref{eq:inter} in terms of the sub-systems~$\dot x(t)=f(t,x)$ and~$\dot x(t)=g(t,x)$.
\end{Remark}

The next result follows
immediately from   Thm.~\ref{thm:rusgens}.

\begin{Corollary}
\label{coro:alpcon}  
Suppose that~\eqref{eq:dsys}
is~$\alpha$ contracting. Then any compact and strongly invariant set has Hausdorff dimension smaller than~$\alpha$. 
\end{Corollary}

\begin{Example}
Consider the system  $\dot{x}(t) =  A(t)x(t)$, with
$
A(t) = \begin{bmatrix}
0&-1&0 \\
1&0&0 \\
0&0&-t
\end{bmatrix}
$.  
Take~$\alpha=2+s$ with $ s\in (0,1)$. By Prop.~\ref{prop:explicitak},
$
A^{[2]}=
\begin{bmatrix}
0&0&0\\
0&-t&-1\\
0&1&-t
\end{bmatrix}
$ and $A^{[3]}= -t$. Hence, the $\alpha$-additive compound is
\[
A^{[\alpha]}=((1-s)A^{[2]}) \oplus (sA^{[3]}) =  \begin{bmatrix}
-st&0&0\\
0&-t&s-1\\
0&1-s&-t
\end{bmatrix}.
\]
Note that $\mu_2(A^{[2]}) = 0$, and $\mu_2(A^{[\alpha]})=-st < 0$ for any~$t>0$. That is, this system is $(2+s)$ contracting with $s \in (0,1)$.
Thm.~\ref{thm:rusgens} thus guarantees that any compact and strongly invariant set~$K$ satisfies~$\dim_H{K}<2+s$. Since $s\in (0,1)$ can be arbitrarily small,   \be\label{eq:gpot}
 \dim_H{K} \leq 2.
 \ee
 For example, the set
$
{K} := \{x \in \R^3: x_1^2+x_2^2 \leq  c,\;x_3 =0   \}
$ with any~$c\geq 0$ is compact, strongly invariant, and satisfies~\eqref{eq:gpot}.
\end{Example}

The next results shows that if the system is~$\alpha$ contracting  w.r.t.  $|\cdot|_p$, for some~$p\in\{1,2,\infty\}$, then 
it is also $\bar \alpha $ contracting w.r.t. the same norm
for any~$\bar \alpha \geq \alpha$.
\begin{Theorem}\label{thm:alpbet}
Consider the system $\eqref{eq:dsys}$. Suppose that condition~\eqref{eq:alphacon} holds for some $\mu_p$ with $p = \{1, 2, \infty\}$, and some $\alpha \in [1, n)$. Then~$\eqref{eq:dsys}$ is   $\beta$~contracting for any $\beta \in (\alpha, n]$. 
\end{Theorem}

\begin{IEEEproof}
Consider first  the case that $\alpha$ is an integer, that is, $\alpha = k \in \{1, \dots, n-1\}$. Then~\eqref{eq:alphacon}  becomes 
\[
\mu_p(J_f^{[k]}(t,x))\leq -\eta<0,\text{ for all } t\geq 0, \; x\in D.
\]
Fix arbitrary~$x \in D$ and~$t\geq 0$. To simplify the notation, we write~$J_f$ for~$J_f(t,x)$. 
  Prop.~\ref{prop:monomeak} ensures that $\mu_p(J_f^{[k+1]}  ) \leq \mu_p(J_f^{[k]} )$. Fix~$\varepsilon \in (0, 1  )$. By Corollary~\ref{coro:ineq}, 
\begin{align*}
 \mu_p(J_f^{[\alpha + \varepsilon]} ) 
=& \mu_p(J_f^{[k + \varepsilon]} ) 
= (1-\varepsilon)\mu_p(J_f^{[k]} ) + \varepsilon \mu_p(J_f^{[k+1]} ) \\
= & \mu_p(J_f^{[k]} )  
 - \varepsilon \left (\mu_p(J_f^{[k]} ) - \mu_p(J_f^{[k+1]}  ) \right)\\
\leq &  \mu_p(J_f^{[k]} ),
\end{align*}
so and the system is~$\alpha+\varepsilon$ contracting. 

Now suppose that $\alpha$ is not an integer, i.e.~$\alpha=k+s$, with~$k$ an integer and~$s\in(0,1)$. Then condition~\eqref{eq:alphacon}  becomes 
\[
(1-s)\mu_p(J_f^{[k]} ) + s \mu_p(J_f^{[k+1]} ) \leq -\eta <0.
\]
We claim that
\be\label{eq:inpart}
  \mu _p (J_f^{[k]}) \geq \mu_p (J_f^{[k+1]}).
  \ee
  Indeed, 
  if $\mu_p(J_f^{[k]} )  \geq 0 $, then $\mu_p(J_f^{[k+1]} ) \leq -\eta/s< 0$, so~\eqref{eq:inpart} holds, 
  and if $\mu_p(J_f^{[k]} ) < 0$, then~\eqref{eq:inpart} follows from Prop.~\ref{prop:monomeak}.  Hence, for any~$\varepsilon \in (0, 1-s )$,
\begin{align*}
 \mu_p(J_f^{[\alpha + \varepsilon]} )  
=& (1- s -\varepsilon)\mu_p(J_f^{[k]} )  
 + (s +\varepsilon) \mu_p(J_f^{[k+1]} ) \\
=& \mu_p(J_f^{[\alpha]} ) 
-\varepsilon \left(\mu_p(J_f^{[k]} ) - \mu_p(J_f^{[k+1]} ) \right) \\
\leq & \mu_p(J_f^{[\alpha]} ),
\end{align*}
and this completes the proof.
\end{IEEEproof}
 
Theorem~\ref{thm:alpbet} implies the following result. 

\begin{Corollary}
Consider the dynamical  system $\eqref{eq:dsys}$. Suppose that condition~\eqref{eq:alphacon} holds for some $\mu_p$ with $p = \{1, 2, \infty\}$, and some $\alpha \in [1, n)$. Then there exists a minimal real value~$\alpha^*\in[1,\alpha]$
such 
that~$\eqref{eq:dsys}$ is   $\beta$~contracting for 
any~$\beta > \alpha^*$. 
\end{Corollary}
In other words, contraction  is not a binary property, but rather   
the system is located  on a continuous axis of contraction  level.  It is important to note that the value~$\alpha^*$
depends on the norm that induces the matrix measure.
This is also true of standard contraction, where the analysis of contraction critically depends on using the ``right'' norm. 

Several recent papers considered systems that are, in some sense, 
on ``the verge of contraction''~\cite{LOHMILLER1998683,transverse2014,sontag2014three,jafarpour2020weak,3gen_cont_automatica,cast_book}. Such systems are referred to as semi-contracting~\cite{LOHMILLER1998683,wensing}, or sometimes weakly-contracting~\cite{jafarpour2020weak} (note that this terminology is used instead for $k$-order contraction in~\cite{weak_slotine}). Since $1$-order contraction corresponds to contracting systems, we can expect semi-contracting systems to be~$\alpha$ contracting for some~$\alpha>1$.
This is indeed the case. We demonstrate this for the important example of studying  synchronization
using contraction theory~\cite{Wang,Pham2007}.
\begin{Example}
Consider the LTI system
\be\label{eq:lap}
\dot x= -Lx,
\ee
  where~$L$ is   the Laplacian
 of a   (directed or undirected) weighted graph with a globally reachable vertex.
We claim that~\eqref{eq:lap} 
is not $1$-order contracting w.r.t.   any norm. 
Yet, for any~$\varepsilon \in(0,1)$  
  there exists a  vector norm~$|\cdot|$ such that~\eqref{eq:lap}
is~$1+\varepsilon$ contracting w.r.t. to~$|\cdot|$. 

Indeed, for any~$c\in\R$ we have that~$c 1_n$ is an equilibrium of~\eqref{eq:lap}, so the system cannot be $1$-order contracting w.r.t. any norm. 
On the other-hand, the eigenvalues $\lambda_i(A)$, ordered as in~\eqref{eq:ord_eig},
satisfy~$\lambda_1=0$
and~$\real(\lambda_2)<0$.
Fix~$\varepsilon\in(0,1)$. By~\eqref{eq:specalpha}, 
\[
\real(  \lambda_1(A^{[1+\varepsilon]}) )
  = \real(   \lambda_1(A) + \varepsilon \lambda_{2}(A))<0,
\]
so~$A^{[1+\varepsilon]}$ is Hurwitz,
and it is well-known~\cite{sontag_cotraction_tutorial,Wang2add} that this implies that there exists a matrix measure~$\mu$ such that~$\mu ( A^{[1+\varepsilon]})<0$.
Combining this with Corollary~\ref{coro:alpcon}   implies that any compact and strongly invariant set~$K$ of the dynamics satisfies~$\dim_H K\leq 1$. This agrees with the fact that  the dynamics converges to ``lines''.
\end{Example}

The next example demonstrates an application of our theoretical results to the control of a chaotic system. 

\begin{figure}[t]
 \begin{center}
\includegraphics[width=8cm,height=6cm]{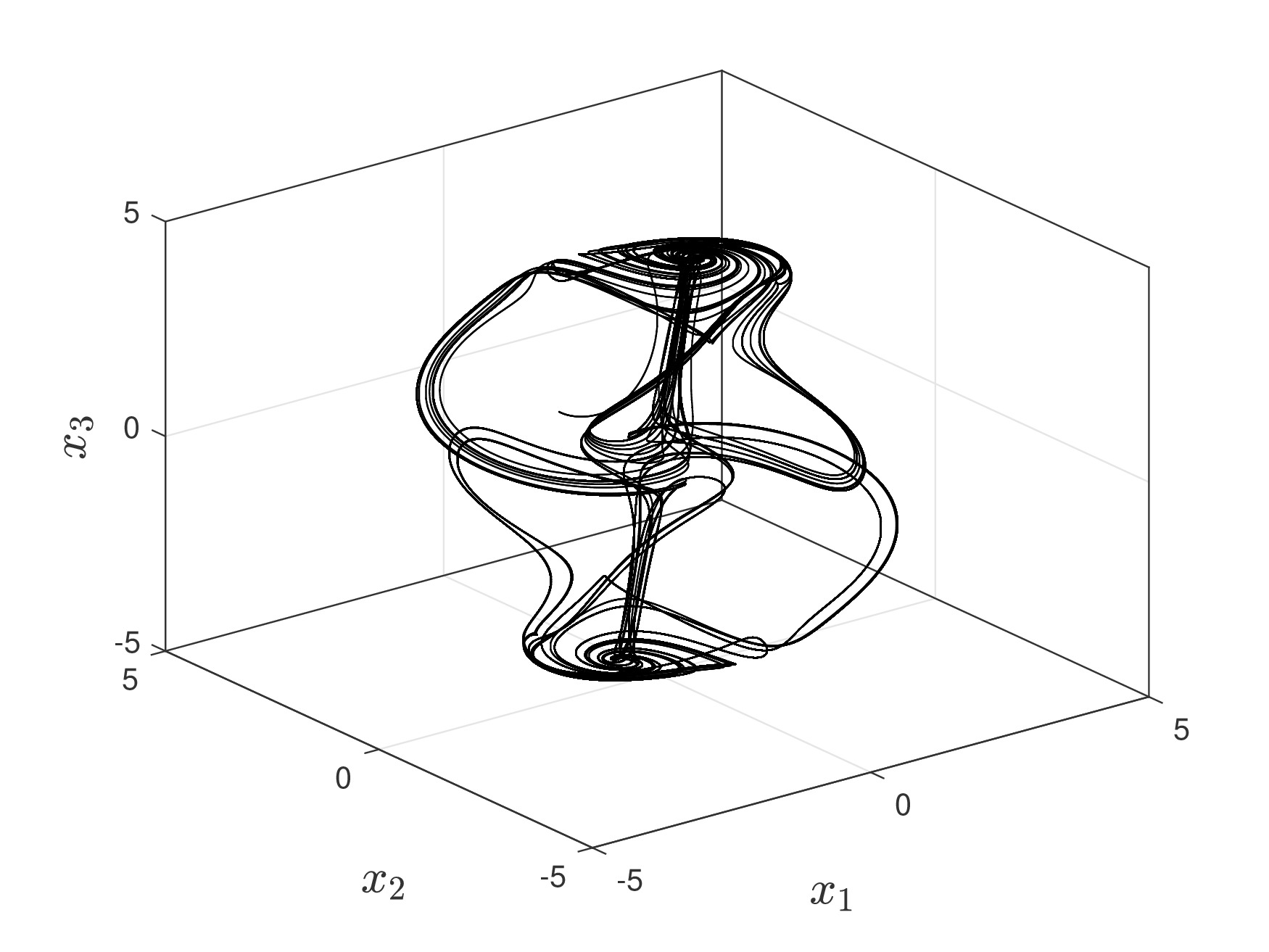}
\caption{A trajectory of~\eqref{eq:thom} emanating from~$x(0)=\begin{bmatrix} 
 -1& 1&1 \end{bmatrix}^T$
for the dissipation parameter in~\eqref{eq:bstar}. }\label{fig:chaos}
\end{center}
\end{figure}

\begin{figure}[t]
 \begin{center}
\includegraphics[width=8cm,height=6cm]{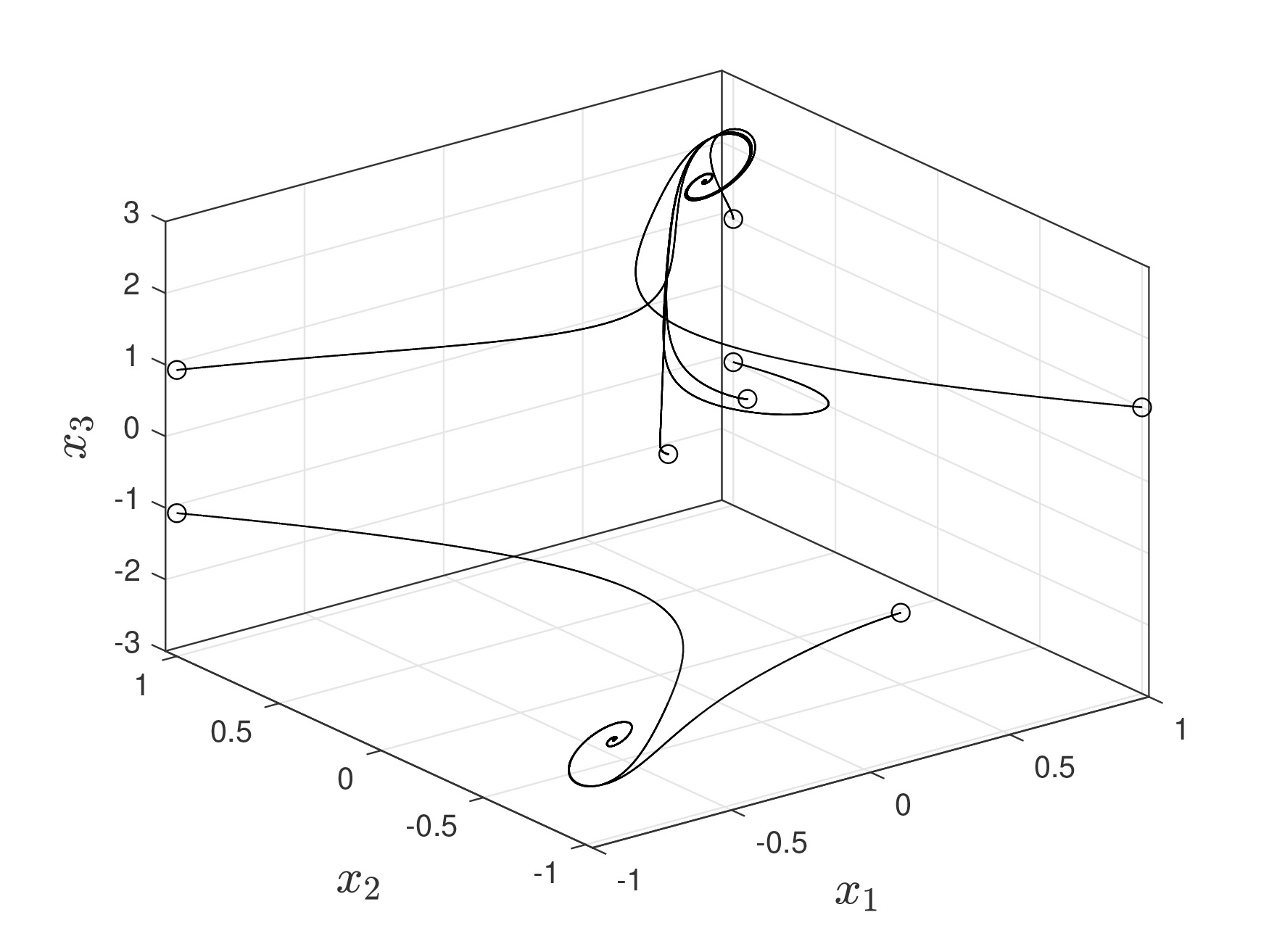}
\caption{Several trajectories of  the closed-loop system for the dissipation parameter   in~\eqref{eq:bstar}. The circles denote the initial conditions of the trajectories.  }\label{fig:chaos_closed}
\end{center}
\end{figure}

\begin{Example}
A popular example for a chaotic system, introduced by Thomas~\cite{thomas99} (see also the recent review~\cite{chaos_survey}),
is 
Thomas' cyclically symmetric attractor:
\begin{align} \label{eq:thom}
\dot x_1 =&  \sin(x_2)-bx_1 , \nonumber \\
\dot x_2 =&  \sin(x_3) - bx_2, \\
\dot x_3 =& \sin(x_1) - bx_3, \nonumber
\end{align}
where $b>0$ is the
dissipation  constant. 
Note that the convex set~$D: = \{x\in\R^3: b |x|_\infty \leq 1 \}$ is an
invariant set of the dynamics.

This system undergoes a series of bifurcations as~$b$ decreases. 
 For $b > 1$   the origin is the single stable equilibrium. When~$ b = 1$
  it undergoes a pitchfork bifurcation, splitting into two attractive fixed points. As~$b$ is decreased further to~$b \approx 0.32899$ these undergo a Hopf bifurcation, creating a stable limit cycle. The limit cycle  undergoes a period doubling cascade and becomes chaotic at~$b  \approx 0.208186$.  
 
 Fig.~\ref{fig:chaos}
 depicts  the solution of the system emanating from~$\begin{bmatrix}-1& 1&1\end{bmatrix}^T$ for
 \be\label{eq:bstar}
 b=0.193186
\ee
 Note the symmetric strange attractor. 
 
 Let~$f$ denote the vector field in~\eqref{eq:thom}.
The Jacobian  is
\begin{align*}
   J_f(x)=\begin{bmatrix}-b&\cos(x_2)&0 \\ 0&-b&\cos(x_3) \\ \cos(x_1)&0&-b\end{bmatrix},
   \end{align*}
   and thus
 \begin{align*}
   J_f^{[2]}(x)=\begin{bmatrix}-2b&\cos(x_3)&0 \\ 0&-2b&\cos(x_2) \\ -\cos(x_1)&0&-2b\end{bmatrix},
   \end{align*}
and~$
   J_f^{[3]}=\trace(J_f(x))=-3b$. 
This implies that the system is $3$ contracting (that is, dissipative), w.r.t. any norm, for any~$b>0$. Let~$\alpha = 2+s$, with~$s\in(0,1)$. Then
\begin{align*}
    &J_f^{[\alpha]}(x)
    =(1-s)J_f^{[2]}(x)\oplus sJ_f^{[3]}(x)\\
    &=
    \begin{bmatrix}
    -(2+s)b & (1-s)\cos(x_3) & 0 \\
    0 & -(2+s)b & (1-s)\cos(x_2) \\
    -(1-s)\cos(x_1) & 0 & -(2+s)b \\
    \end{bmatrix}.
\end{align*}
This implies that
\be\nonumber
\mu_{1 }(J_f^{[\alpha]}(x)) \leq 1-2b-s(b+1), \text{ for all } x\in D.
\ee
We conclude that for any~$b\in(0,1/2) $
 the system is~$2+s $ contracting for any
$
    s>\frac{1-2b}{1+b}
$.

We now demonstrate how our results can be applied to design a partial-state   controller for the system guaranteeing that the closed-loop system has a ``well-ordered'' behaviour. 
Suppose that 
the closed-loop system is:
\[
\dot x = f(x) + g(x),
\]
where~$g $ is the controller.
Let~$\alpha = 2 + s$, with~$s \in(0,1)$. 
The Jacobian of the closed-loop system is~$J_{cl}:=J_f+J_g$,  so
\begin{align*}
    \mu_1(J_{cl}^{[\alpha]})& = \mu_1(J_f^{[\alpha]}+J_g^{[\alpha]})  \leq \mu_1(J_f^{[\alpha]})+\mu_1(J_g^{[\alpha]}) \\
    &\leq 1-2b-s(b+1)+\mu_1(J_g^{[\alpha]}).
\end{align*}
This implies that the closed-loop system is~$\alpha$ contracting if
\begin{align}\label{eq:cond_cont}
    \mu_1(J^{[\alpha]}_g (x) ) <    s(b+1)+2b-1
    \text{ for all } x\in D . 
\end{align}
Consider, for example,   the controller  
\[
g(x_1,x_2)=\diag(c,c,0) x ,\text{ with } c<0.
\]
Then
\[
J_g^{[\alpha]} = c\diag(  2  ,1+ s,1+ s )
\]
and for any~$c<0$  condition~\eqref{eq:cond_cont} becomes
\begin{align}\label{eq:cond_cont1}
     (1+ s)c <    s(b+1)+2b-1.
\end{align}
This provides a simple recipe  for 
determining the gain~$c$ so that the closed-loop system is~$2+s$ contracting. For example,  
  when~$s \to 0$, Eq.~\eqref{eq:cond_cont1} yields
\[
c<2b-1
\]
and this guarantees that the closed-loop system is~$2$-order contracting. Recall that in a~$2$-order contracting system 
every nonempty omega limit set is a single equilibrium,
 thus ruling out chaotic attractors and even non-trivial limit cycles~\cite{li1995}.  
Fig.~\ref{fig:chaos_closed} depicts the behaviour of the closed-loop system with~$b$ as in~\eqref{eq:bstar} and~$c=2b-1.1$. The closed-loop system is thus $2$-order contracting, and as expected  
every solution converges to an equilibrium. 
\end{Example}

 \section{Conclusion}
 The~$k$~multiplicative and~$k$~additive compounds of a matrix play an important role in geometry, multi-linear algebra, dynamical systems, and more. These compounds are based on~$k\times k$ minors and are thus  defined for integer values of~$k$ only. The $k$ compounds were recently used to study an extension of contracting systems to~$k$-order contracting  systems~\cite{kordercont}. 
 
 Here, we generalised $k$ compounds to~$\alpha$   compounds, with~$\alpha$ real, and analyzed the properties of these compounds.
 As an application, we showed that these compounds provide more direct and  intuitive interpretation of important functions, e.g.~$\omega(K, \alpha, g)$ appearing in the seminal work of    Douady 
 and Oesterl\'{e}. We also introduced the new   notion of~$\alpha$ contracting systems, with~$\alpha$ real, generalizing the  
 notion of~$k$-order contracting systems with~$k$ an integer, recently analyzed in~\cite{kordercont}. Thus, rather than a binary choice $-$
 contracting or not contracting in a given metric $-$ one can place any 
 system on a continuous axis of contraction levels.

  Due to space limitations, we focused here on theoretical results, but  we believe that  many applications are possible. 
  First, there exist  nonlinear systems where the ``level of contraction'' naturally changes in a continuous way, for example, systems that involve a continuous-time dynamics and discrete-time switching (see, e.g.~\cite{doi:10.1002/aic.690460317}). 
Second, contraction theory (i.e. the theory of 1-order contracting systems~\cite{LOHMILLER1998683}) has found many applications in control synthesis (see e.g., \cite{sanfelice2011convergence,slotine2005study,manchester2018unifying,russo2009solving,RFM_entrain,entrain2011, wu2019robust}). An interesting research direction is to apply the  generalization described here to 
control synthesis in such contexts. Finally, our results could be used to define generalized notions of convexity in optimization and machine learning. Just as Riemannian convexity of a scalar function with respect to a metric is equivalent to contraction in that metric of natural gradient descent~\cite{wensing}, notions of $\alpha$ Riemannian convexity could similarly be defined through equivalent $\alpha$ contracting autonomous dynamical systems.

 \end{document}